\documentclass[11pt]{article}
\usepackage{amsmath,amsthm,amssymb,mathrsfs}
\usepackage{bm,pifont,xcolor, tikz}

\usepackage[normalem]{ulem}

\makeatletter
\newcommand{\rmnum}[1]{\romannumeral #1}
\newcommand{\Rmnum}[1]{\expandafter\@slowromancap\romannumeral #1@}
\makeatother

\usepackage{geometry}
\usepackage[colorlinks=true,
linkcolor=blue,citecolor=blue,
urlcolor=blue]{hyperref}

\makeatletter
\def\@seccntDot{.}
\def\@seccntformat#1{\csname the#1\endcsname\@seccntDot\hskip 0.5em}
\renewcommand\section{\@startsection{section}{1}{\z@}%
{18\p@ \@plus 6\p@ \@minus 3\p@}%
{9\p@ \@plus 6\p@ \@minus 3\p@}%
{\large\bfseries\boldmath}}
\renewcommand\subsection{\@startsection{subsection}{2}{\z@}%
{12\p@ \@plus 6\p@ \@minus 3\p@}%
{3\p@ \@plus 6\p@ \@minus 3\p@}%
{\bfseries\boldmath}}
\renewcommand\subsubsection{\@startsection{subsubsection}{3}{\z@}%
{12\p@ \@plus 6\p@ \@minus 3\p@}%
{\p@}%
{\bfseries\boldmath}}
\makeatother

\usepackage{microtype}

\theoremstyle{plain}
\newtheorem{theorem}{Theorem}[section]
\newtheorem{lemma}{Lemma}[section]

\newtheorem{proposition}{Proposition}[section]

\theoremstyle{definition}

\numberwithin{equation}{section}
\allowdisplaybreaks
\title{Uniquely $K_{4}$-saturated graphs with large maximum degree
}

\usepackage{authblk}

\author[1\thanks{ Email: yuyingli1205@163.com}]
{Yuying Li}

\author[2,3\thanks{ Email: weiyuanzhang@nuaa.edu.cn}]
{Weiyuan Zhang}

\author[2,3\thanks{Corresponding author. Email: kexxu1221@126.com}]
{Kexiang Xu}

\affil[1]{\footnotesize School of Mathematical Science, Jiangsu  Second Normal University,
\protect\\ Nanjing, Jiangsu, 210013, PR China}

\affil[2]{\footnotesize School of Mathematics,
	Nanjing University of Aeronautics and Astronautics,
\protect\\ Nanjing, Jiangsu, 210016, PR China}

\affil[3]
{\footnotesize
MIIT Key Laboratory of Mathematical Modelling and High Performance Computing of Air  Vehicles, Nanjing, Jiangsu, 210016, PR China
}

\date{}

\begin{document}
\maketitle

\begin{abstract}
A graph $G$ is  \textit{uniquely $H$-saturated} if it contains no copy of $H$ as a subgraph, but adding any edge to $G$ creates exactly one copy of $H$ in the resultant graph.
In this paper,   we consider the uniquely  $K_{4}$-saturated graphs and characterize the uniquely  $K_{4}$-saturated graphs on $n$  vertices  with maximum degree $\Delta\geq n-5$.
We prove that if $G$  is a uniquely  $K_{4}$-saturated graph with universal vertex $x$, then $G-x$ is either a star or a Moore graph of diameter 2.
The uniquely  $K_{4}$-saturated graphs with maximum degrees $n-3$ and $n-5$ are $\overline{C_{7}}$ and $G_{10}$  in  Figure 2,  respectively.
But there are no uniquely  $K_{4}$-saturated graphs with maximum degree $n-2$ or $n-4$.

\par\vspace{2mm}

\par\vspace{1mm}

\end{abstract}

\section{Introduction}

All graphs considered in this paper are finite,  undirected and simple.  Let $G=(V(G), E(G))$  be a graph,  as usual,  denote by $e(G)$ and $\overline{G}$  the number of edges and the complement of $G$, respectively.
For any $v\in V(G)$, let $d_{G}(v)$  and $N_{G}(v)$  denote its degree and neighborhood  in $G$, respectively.
Let $N_{G}[v]=N_{G}(v)\cup\{v\}$  be the  closed  neighborhood of $v$ in $G$.
A vertex $v$ is a  \textit{universal vertex} of  $G$ if  $d_{G}(v)=|V(G)|-1$,  and a vertex $v$ is an  \textit{isolated vertex} of  $G$ if  $N_{G}(v)=\emptyset$.
As usual, let $\Delta(G)$  be the maximum degree of graph $G$.  The  \textit{distance $d_{G}(x,y)$}  of two vertices $x$, $y$ is the length of a shortest $x-y$ path in $G$.  The greatest distance between any two vertices in $G$ is the \textit{diameter} of $G$, denoted by $diam(G)$.
The \textit{girth} $g(G)$ of a graph $G$ is the minimum length of  cycles in it.  The graphs with diameter $d$ and girth $2d+1$ are the \textit{Moore graphs}.
 Denote by $K_{n}$, $P_{n}$ and $C_{n}$ the complete graph, path and cycle on $n$ vertices, respectively.
For $A\subseteq V(G)$, let $G[A]$   be the subgraph of $G$ induced by $A$.  Given vertex sets $A, B \subseteq V(G)$, we denote by $E_{G}(A,B)$, or $E(A,B)$ for simplicity, the  edge set between $A$ and $B$ in $G$  and let $A\setminus B$ denote $A-B$.  For any edge $e\in E(\overline{G})$ (here we call $e$ a non-edge of $G$),  we  denote by $G+e$ the graph obtained from $G$ by adding the new edge $e$.  Given any two  vertex-disjoint graphs $G$ and $H$,  their  \textit{union}  $G\cup H$   is the graph with vertex set $V(G)\cup V(H)$ and edge set $E(G)\cup E(H)$  and  their  \textit{join}  $G\vee H$   is the graph with vertex set $V(G)\cup V(H)$ and edge set $E(G)\cup E(H)\cup\{gh : g\in V(G), h\in V(H)\}$.  For brevity, we write $[t]=\{1,2,\ldots,t\}$  for a positive integer  $t$.  Please refer to  \cite{A.Bondy2008}  for  more definitions and  notations.

Given a graph $H$,  a graph $G$ is \textit{$H$-free} if $G$ does not contain $H$ as a  subgraph.
A  graph $G$  is   \textit{$H$-saturated} if $G$ is $H$-free, but $G+e$  contains a copy of $H$ for any $e\in E(\overline{G})$.
The study of $H$-saturated graphs began when Mantel \cite{W.Mantel1907} proved that the $K_{3}$-saturated graph on  $n$  vertices with the most edges is $K_{\lceil\frac{n}{2}\rceil,\lfloor\frac{n}{2}\rfloor}$.
Tur\'{a}n \cite{P.Turan1941} generalized this  result, proving that the  $K_{r}$-saturated graphs on  $n$ vertices  with the most edges are the complete  $(r-1)$-partite graphs with partite sets as balanced as possible.
In the opposite direction,  Erd\H{o}s \textit{et al.} \cite{P.Erdos1964} determined that the $K_{r}$-saturated graphs on  $n$ vertices  with the fewest edges   are   $K_{r-2}\vee\overline{K_{n-r+2}}$.
Lei \textit{et al.} \cite{LOSWZ2019} established sharp lower and upper bounds on the number of edges in $n$-vertex $\mathcal F$-saturated graphs, where $\mathcal F$ is the family of $k$-edge-connected graphs.
Recently Shi \textit{et al.} \cite{Huang2024, Lan2025} provided the bounds on the minimum number of edges in $K_{3,3}$-saturated graphs  and  $C_{6}$-saturated graphs, respectively.
Please see  an informative survey  \cite{B.L.Currie2021}  for some detailed results about the $H$-saturated graphs  with the fewest edges.

For a graph $H$ and an $H$-saturated graph $G$,  we say that $G$ is \textit{uniquely $H$-saturated} if  $G+e$  contains exactly one copy of $H$ for any $e\in E(\overline{G})$.
Every complete graph with fewer than $|V(H)|$  vertices is trivially  uniquely $H$-saturated, since there are no edges to be added to it.
A uniquely $H$-saturated graph is \textit{nontrivial} if it has at least  $|V(H)|$ vertices.
The graphs found in \cite{P.Erdos1964} are uniquely $K_{r}$-saturated,  but the graphs  with at least five vertices found in \cite{P.Turan1941} are  not.
Questions about the uniquely $H$-saturated graphs focus on their existence.
Cooper  \textit{et al.} \cite{J.Cooper2012}  initiated the study of the uniquely saturated graphs by determining all uniquely $C_{k}$-saturated graphs  for  $k\in\{3, 4\}$.
They proved that there are exactly ten uniquely $C_{4}$-saturated graphs   and  a graph is uniquely $C_{3}$-saturated if and only if it is a star or a Moore graph of diameter 2.
Later, Wenger and West \cite{P.S.Wenger2016}  proved that the  nontrivial  uniquely $C_{5}$-saturated graphs are precisely the friendship graphs and there are no  nontrivial  uniquely $C_{k}$-saturated graphs for $k\in\{6,7,8\}$.
Furthermore, they showed that there are only finitely many  uniquely $C_{k}$-saturated graphs for all $k\geq6$  and conjectured that no such graphs exist.
Recently Li  \textit{et al.} \cite{Li2026} characterized  the  uniquely $C_{4}^{+}$-saturated graphs based on girth and the number of triangles, where $C_{4}^{+}$  is the diamond graph consisting of adding an edge in  cycle $C_{4}$.
Ji  \textit{et al.} \cite{Ji2026} proved that there are no nontrivial uniquely bowtie-saturated graphs.
For the path, Wenger \cite{P.S.Wenger2010}  proved  that  no nontrivial uniquely $P_{k}$-saturated graphs exist  for $k\geq5$  and also characterized the nontrivial uniquely $P_{k}$-saturated graphs  for  $k\in\{2,3,4\}$.
Some results on uniquely $T$-saturated graphs can be  referred to \cite{L.W.Berman2016} for tree $T$.
We remark that by now  few examples of uniquely $K_{r}$-saturated graphs are known, and little is known about their properties.
Hartke  and  Stolee  \cite{S.G.Hartke2012} considered  the  uniquely $K_{r}$-saturated graphs  with  no universal vertex.  They  found ten new uniquely $K_{r}$-saturated graphs with $4\leq r\leq7$  by  computational methods, as well as two new infinite families of uniquely $K_{r}$-saturated  Cayley graphs.
Some  results about uniquely saturated graphs on the complete uniform hypergraph can be found in \cite{A.Gyarfas2018}.

So far,  the uniquely $C_{4}$-saturated graphs and the uniquely $C_{4}^{+}$-saturated graphs have been characterized.
But  the research  on the uniquely $K_{4}$-saturated graphs is quite challenging (most of the uniquely $K_4$-saturated graphs were discovered using computational methods).
In this paper,  we  consider   the uniquely  $K_{4}$-saturated graphs with large maximum degree.
In Section 2,   we prove  some  structural properties  of the uniquely $K_{4}$-saturated graphs  and  determine the uniquely  $K_{4}$-saturated graphs with  universal  vertex.
In Section 3,  we show  that there are no uniquely  $K_{4}$-saturated graphs with maximum degree $n-2$ or $n-4$  and   characterize  the uniquely  $K_{4}$-saturated graphs with maximum degrees $n-3$ and $n-5$.

\section{Structural properties of the uniquely $K_{4}$-saturated graphs}

For  complete graph  $K_{r}$,    only the uniquely $K_{3}$-saturated graphs are  fully characterized  in   \cite{J.Cooper2012}.

\begin{theorem}[\cite{J.Cooper2012}]\label{theorem:2.1}
 $G$ is a uniquely $K_{3}$-saturated graph if  and  only  if  it is  a star or a Moore graph of diameter 2.
\end{theorem}

\begin{proposition}\label{proposition:2.1}
Let  $G$  be a uniquely $K_{r}$-saturated graph.

(\rmnum{1}).  If  $x$  is a universal  vertex  of  $G$,  then  $G-x$  is a uniquely $K_{r-1}$-saturated graph.

(\rmnum{2}).   $G\vee K_{1}$  is a uniquely $K_{r+1}$-saturated graph.
\end{proposition}
\begin{proof}
(\rmnum{1}).  Since $x$  is a universal  vertex  of  $G$,  for  any  non-edge $e$  of  $G-x$,    $e$  is also a non-edge of  $G$  and  $x$  belongs to  the  only $K_{r}$  in  $G+e$.  Thus  $(G-x)+e$   contains  exactly one $K_{r-1}$,   which   means that  $G-x$  is a uniquely $K_{r-1}$-saturated graph.

(\rmnum{2}).  For  any  non-edge $e$  of  $G\vee K_{1}$,    $e$  is also a non-edge of  $G$.   Then  $(G\vee K_{1})+e$     contains  exactly one $K_{r+1}$  since    $G+e$   contains  exactly one $K_{r}$  as a subgraph.
\end{proof}

By  \autoref{theorem:2.1}  and  \autoref{proposition:2.1},  we  can get the following result about  the  uniquely $K_{4}$-saturated graphs  with   universal  vertex.

\begin{theorem}\label{theorem:2.2}
Let  $G$  be a  graph  with a  universal  vertex  $x$, then  $G$ is a uniquely $K_{4}$-saturated graph  if  and  only  if  $G-x$   is  a star or a Moore graph of diameter 2.
\end{theorem}

A uniquely $K_{r}$-saturated graph  with no  universal  vertex   denotes  as  \textit{$r$-primitive  graph}  by  \cite{S.G.Hartke2012}.   Collins and  Kay discovered the only previously known infinite family of  $r$-primitive graphs   $\overline{C_{2r-1}}$.  Collins and Cooper discovered two more $4$-primitive graphs of orders $10$ and $12$,  respectively.  In \cite{S.G.Hartke2012},  using a computational method,   Hartke and Stolee  exhaustively searched the  uniquely $K_{4}$-saturated graphs of order at most $20$  and   found four new $4$-primitive graphs  of orders $13$ and $18$,  respectively. By  now only $7$ $4$-primitive graphs  are  known.

Let the \textit{wheel graph} $W_{n}$ denote the graph  obtained from $C_{n}$ by adding a new vertex and joining it to every vertex of $C_{n}$, that is , $W_{n}=K_{1}\vee C_{n}$. And let $F^{-}(2,4,2)$ (see Figure 1 right) be the graph obtained by deleting the common edge of two copies of $K_{4}$ intersecting in a common edge.
Let us show that the uniquely $K_{4}$-saturated graph  has some  forbidden  subgraphs  in  the  following result.

\begin{center}
\begin{tikzpicture} [inner sep=0.7mm, scale=0.8]
\tikzstyle{place}=[circle,draw=black!,fill=black!]
\tikzstyle{r}=[circle,draw=red!,fill=red!]

\node at (0,0)    [place] {};
\node at (3,0)    [place] {};
\node at (0,3)    [place] {};
\node at (3,3)    [place] {};
\node at (1.5,1.5)    [place] {};

\draw[black,line width=1] (1.5,1.5) -- (0,0);
\draw[black,line width=1] (1.5,1.5) -- (3,0);
\draw[black, line width=1] (1.5,1.5) -- (0,3);
\draw[black, line width=1] (1.5,1.5) -- (3,3);
\draw[black, line width=1] (0,0) -- (3,0);
\draw[black, line width=1] (0,0) -- (0,3);
\draw[black, line width=1] (3,3) -- (3,0);
\draw[black, line width=1] (3,3) -- (0,3);

\node at (6,1.5)    [place] {};
\node at (9,1.5)    [place] {};
\node at (7.5,0.5)    [place] {};
\node at (7.5,2.5)    [place] {};
\node at (7.5,3.5)    [place] {};
\node at (7.5,-0.5)    [place] {};

\draw[black,line width=1] (7.5,2.5) -- (7.5,3.5);
\draw[black,line width=1] (7.5,2.5) -- (6,1.5);
\draw[black, line width=1] (7.5,2.5) -- (9,1.5);
\draw[black, line width=1] (7.5,0.5) -- (7.5,-0.5);
\draw[black, line width=1] (7.5,0.5) -- (6,1.5);
\draw[black, line width=1] (7.5,0.5) -- (9,1.5);

\draw[black,line width=1] (7.5,3.5) -- (9,1.5);
\draw[black,line width=1] (7.5,3.5) -- (6,1.5);
\draw[black, line width=1] (7.5,-0.5) -- (9,1.5);
\draw[black, line width=1] (7.5,-0.5) -- (6,1.5);

\end{tikzpicture}

\end{center}
\begin{center}
Figure 1 : Wheel graph $W_{4}$  (left)  and  Graph $F^{-}(2,4,2)$ (right).
\end{center}

\begin{lemma}\label{lemma:2.1}
Let $G$  be a  uniquely $K_{4}$-saturated graph,  then $G$  contains no  $W_{4}$ and  $F^{-}(2,4,2)$  as  subgraphs.
\end{lemma}
 \begin{proof}
If $G$ contains a copy of $W_{4}$ as a subgraph,  let  $u$ and  $v$ be two nonadjacent  vertices with degree 3  in the $W_{4}$,  then $G+uv$  contains  two $K_{4}$, which  contradicts  the  fact that $G$  is  uniquely $K_{4}$-saturated.
Similarly, if $G$ contains a copy of $F^{-}(2,4,2)$  as a subgraph,  we  consider  the  two vertices $u$,  $v$  with degree 4 in the  $F^{-}(2,4,2)$   and we can get a contradiction again.
 \end{proof}

\begin{lemma}\label{lemma:2.2}
For a $4$-primitive graph $G$,  let $x$  be any vertex  in $G$  and  $A=V(G)\setminus N_{G}[x]$,  then  $G[A]$  has  no  isolated vertex.
\end{lemma}
\begin{proof}
Let  $N_{G}(x)=\{x_{1},\ldots,x_{t}\}$  and  $A=\{y_{1},\ldots,y_{p}\}$.
If  there exists an isolated vertex,  say  $y_{1}$,  in $G[A]$,  then  $y_{1}$  is  not  adjacent to $x$  and $y_{i}$  with  $2\leq i\leq p$.
Now we claim $N_{G}(y_{1})=N_{G}(x)$.  If  not,  we assume that $x_{1}$  is  not  adjacent to  $y_{1}$.  Then  $G+x_{1}y_{1}$  contains one $K_{4}$  with vertex set $\{x_{1}, y_{1}, a, b\}$.
Moreover,  $\{a, b\}\subset N_{G}(x)$  since $y_{1}$   is  not  adjacent to any vertex in $A\cup\{x\}$.  But  $G[\{x, x_{1}, a, b\}]=K_{4}$,  a contradiction.   Thus  $N_{G}(y_{1})=N_{G}(x)$.
 Since  $G+xy_{1}$  contains  exactly one $K_{4}$,  we   have  $G[N_{G}(x)]=P_{2}\cup(t-2)K_{1}$  with   $P_{2}=uv$  and  $t\geq3$.
If  $|A|=1$,    $G+uw$  contains  no $K_{4}$  with  $w\in N_{G}(x)\setminus\{u,v\}$,  a contradiction.  If  $|A|\geq2$,  for  any  $y_{i}\in A$,   $G+xy_{i}$  contains one $K_{4}$  with  vertex set  $\{x, y_{i}, u, v\}$.  Then  $G[A]$  has   no  edge  since  $G$  is  $K_{4}$-free. But  $G+uw$  contains  no $K_{4}$  with  $w\in N_{G}(x)\setminus\{u,v\}$,  a contradiction  again.  Therefore,  $G[A]$  has  no  isolated vertex.
\end{proof}

\section{The uniquely $K_{4}$-saturated graphs with large maximum degree}

In  this  section,  we  consider the uniquely $K_{4}$-saturated graphs  on $n$  vertices  with maximum degree $n-5\leq\Delta\leq n-2$.  We  prove that   there are no uniquely $K_{4}$-saturated graphs  on  $n$  vertices  with  maximum degree  $n-2$  or  $n-4$,   and  characterize the uniquely  $K_{4}$-saturated graphs with maximum degrees $n-3$ and $n-5$.

\begin{theorem}\label{theorem:3.1}
 There are no uniquely $K_{4}$-saturated graphs  on  $n$  vertices  with  maximum degree  $n-2$  or  $n-4$.
\end{theorem}
\begin{proof}
Let  $G$  be a  uniquely $K_{4}$-saturated graph  on  $n$  vertices  with  maximum degree  $\Delta$,  $x$  be a vertex  with  degree  $\Delta$  and $A=V(G)\setminus N_{G}[x]$.   $G[N_{G}(x)]$  is $K_3$-free since $G$ is $K_4$-free.  If  $\Delta=n-2$,  then  $G[A]$  has an isolated vertex,  which contradicts    \autoref{lemma:2.2}.
Thus  there are no uniquely $K_{4}$-saturated graphs  on  $n$  vertices  with  maximum degree  $n-2$.

Now we assume   $\Delta=n-4$,  let  $A=\{u,v,w\}$  and  $N_{G}(x)=\{x_{1},x_{2},\ldots,x_{n-4}\}$.   $G[N_{G}(x)]$  is $K_3$-free since  $G$  is  $K_{4}$-free.   By  \autoref{lemma:2.2}, we have $2\leq e(G[A])\leq3$.

\textbf{Case 1.} $e(G[A])=2$.

Without loss of generality,  we  say that $v$  is   adjacent to $u$  and  $w$.
Since  $G$  is   uniquely $K_{4}$-saturated  and  $x$  is not adjacent to $u$  and  $w$,  $G+xv$  contains exactly one  $K_{4}$   and   $u,w\notin K_{4}$.   Without loss of generality, let  $K_{4}$  in  $G+xv$  have  vertex set $V(K_{4})=\{v,x,x_{1},x_{2}\}$.
Moreover,  $w$ and $u$  have  at most one  adjacent vertex  in $\{x_{1}, x_{2}\}$  since $G$  is  $K_{4}$-free.

If $N_{G}(w)\cap \{x_1, x_2\}=\emptyset$,  then  $G+wx_{1}$  contains one $K_{4}$.   We have   $v\in K_{4}$  and  $u, x_{2}\notin K_{4}$   in $G+wx_{1}$  since  $G[N_{G}(x)]$  is  $K_{3}$-free  and  $u$, $x_{2}$  are  not adjacent to  $w$.   Then   $K_{4}$  in  $G+wx_{1}$  has  vertex set  $\{w,v,x_{1},x_{i}\}$  with  $i\neq2$.  But then $G$  contains  a $W_{4}$ with vertex set $\{v,x,x_1,x_2,x_i\}$,  contradicting \autoref{lemma:2.1}.
Similarly,  $N_{G}(u)\cap \{x_1, x_2\}=\emptyset$ also does not hold.

If $w$ and $u$ have exactly  one adjacent vertex  in $\{x_{1}, x_{2}\}$, respectively,
 then, without loss of generality,  we  say   $N_{G}(w)\cap \{x_1, x_2\}=\{x_{2}\}$.
Now we claim that $N_{G}(w)=\{x_2,x_3,\ldots,x_{n-4},v\}=\big(N_{G}(x)-\{x_{1}\}\big)\cup\{v\}$.
If not,   there  exists $x_{i}$,  $i\geq3$, such  that $x_{i}$  is  not  adjacent to  $w$,  then  $G+wx_{i}$  contains one  $K_{4}$.  We  can  get  $v\in K_{4}$  and  $u\notin K_{4}$ in  $G+wx_{i}$   since  $G[N_{G}(x)]$  is $K_{3}$-free  and  $u$  is  not adjacent to  $w$.   Then   $K_{4}$  in  $G+wx_{i}$  has  vertex set  $\{w,v,x_{i},x_{j}\}$  with  $j\neq1, i$.  But  $G+xv$  contains   two $K_{4}$  with  vertex sets  $\{v,x,x_{1},x_{2}\}$    and  $\{v,x,x_{i},x_{j}\}$,    a contradiction.   Therefore,  we  have $N_{G}(w)=\big(N_{G}(x)-\{x_{1}\}\big)\cup\{v\}$.  Similarly  as above, we can get  $N_{G}(u)=\big(N_{G}(x)-\{x_{i}\}\big)\cup\{v\}$,  $i=1$ or $2$.

If $N_{G}(w)=N_{G}(u)=\big(N_{G}(x)-\{x_{1}\}\big)\cup\{v\}$,  then we  claim that $G[N_{G}(x)\setminus\{x_{1}\}]$  has  exactly  one edge $x_{i}x_{j}$  since  $G+xu$   contains  exactly  one $K_{4}$ with $V(K_4)=\{x,u,x_{i},x_{j}\}$.
Thus  $K_{4}$  in  $G+xu$  and  $G+xw$  have  vertex  sets $\{x,u,x_{i},x_{j}\}$  and  $\{x,w,x_{i},x_{j}\}$,  respectively.
 But   $G+uw$  contains  two $K_{4}$  with  vertex  sets $\{u,v,w,x_{2}\}$  and  $\{u,w,x_{i},x_{j}\}$,  a contradiction.
Thus  $N_{G}(u)=\big(N_{G}(x)-\{x_{2}\}\big)\cup\{v\}$ and $N_{G}(w)=\big(N_{G}(x)-\{x_{1}\}\big)\cup\{v\}$.
If  $G[N_{G}(x)\setminus\{x_{1},x_{2}\}]$  has  an edge $x_{i}x_{j}$,    we have  $G[N_{G}(x)]=2P_{2}\cup(n-8)K_{1}$  with  $2P_{2}=x_{1}x_{2}\cup x_{i}x_{j}$.  Otherwise,  $G+xu$  or  $G+xw$  contains  more than one $K_{4}$,  a contradiction.
 Then  $v$ must belong  to $K_{4}$  in  $G+x_{1}x_{i}$   and  $G+x_{2}x_{j}$  by the  structure of $G$,   which  means that  $v$ is adjacent to  $x_{i}$  and  $x_{j}$.  But  $G$   contains a $K_{4}$  with  vertex set $\{u,v,x_{i},x_{j}\}$, contradicting the fact that  $G$  is  $K_{4}$-free.  Therefore,  $G[N_{G}(x)\setminus\{x_{1},x_{2}\}]$  has  no edge.
Moreover, we can get  that $K_{4}$  in   $G+xu$  has vertex  set $\{x,u,x_{1},x_{i}\}$   and  $K_{4}$  in   $G+xw$  has vertex  set $\{x,w,x_{2},x_{j}\}$  with  $i\neq 2$, $j\neq 1$  and $i\neq j$.  For convenience,  we  say  $i=3$  and  $j=4$.   Then  $G[N_{G}(x)]=P_{4}\cup(n-8)K_{1}$  with  $P_{4}=x_{3}x_{1}x_{2}x_{4}$   since $G$ is $K_4$-free.  Moreover,   $v$ is not adjacent to $x_3, x_4$ since $G+xv$   contains  exactly one $K_4$ with $V(K_4)=\{x,v,x_1,x_2\}$.  But then $G+x_{3}x_{4}$   contains  no  $K_{4}$, a contradiction.

\textbf{Case 2.} $e(G[A])=3$.

In this  case we have  $G[A]=K_{3}$,  then  $u$, $v$ and  $w$  have  no  common adjacent vertex in  $N_{G}(x)$  since  $G$ is  $K_{4}$-free.  Without loss of generality,  let  $K_{4}$  in  $G+xv$  have  vertex set $V(K_{4})=\{v,x,x_{1},x_{2}\}$.   Since  $G$  is  $K_{4}$-free,  $K_{4}$  in  $G+xu$ and $K_4$ in $G+xw$  contain at most one vertex in $\{x_{1}, x_{2}\}$,  respectively.

If  $K_{4}$  in  $G+xu$ and $K_{4}$  in  $G+xw$  contain exactly  one vertex in $\{x_{1}, x_{2}\}$,  then, without loss of generality,   we can  assume  $V(K_{4})=\{x,u,x_{1},x_{4}\}$  in  $G+xu$  and  $V(K_{4})=\{x,w,x_{2},x_{3}\}$  in  $G+xw$  since $G$  is  $K_{4}$-free.
But $G$  contains  a $W_{4}$ with vertex set $\{u,v,w,x_1,x_2\}$,  contradicting \autoref{lemma:2.1}.

If  $K_{4}$  in   $G+xw$  contains exactly  one vertex in $\{x_{1}, x_{2}\}$  and  $K_{4}$  in   $G+xu$  contains no vertex in $\{x_{1}, x_{2}\}$,  then, without loss of generality,   we say  $V(K_{4})=\{x,w,x_{2},x_{3}\}$  in  $G+xw$.  Moreover, we claim that $u$ is not adjacent to $x_3$.  If not, $G$  contains a $W_{4}$ with vertex set $\{u,v,w,x_2,x_3\}$,  contradicting \autoref{lemma:2.1}.
Thus, without loss of generality,   we say  $V(K_{4})=\{x,u,x_{4},x_{5}\}$  in  $G+xu$.
Since  $G$  is  $K_{4}$-free,  $v$  and  $w$  have at most one  adjacent vertex  in $\{x_{4}, x_{5}\}$.
If  $v$  is not adjacent to  $x_{4}$  and  $x_{5}$,  then  $G+vx_{4}$  contains one $K_{4}$  and let $V(K_4)=\{v,x_4,a,b\}$.
 Moreover,  we  claim that $w\in K_{4}$  in   $G+vx_{4}$.  If not, we have $\{a,b\}\subseteq N_{G}(x)$  or $u\in K_{4}$  in   $G+vx_{4}$,  then $G$ contains $K_4$ with vertex set $\{x,x_4,a,b\}$  or  $G+ux$  contains  two  $K_{4}$,  a contradiction.  Thus $w\in K_{4}$ in   $G+vx_{4}$,  which means that $w$  is adjacent to  $x_{4}$.  Similarly,   we can  get  that  $G+vx_{5}$  contains one $K_{4}$  and   $w\in K_{4}$,  which  means that  $w$  is adjacent to  $x_{5}$.  But  $G$  contains a  $K_{4}$  with vertex  set  $\{u,w,x_{4},x_{5}\}$,  a contradiction.  Therefore,  $v$  has exactly   one  adjacent vertex  in $\{x_{4}, x_{5}\}$.   Similarly  as above proof, we can also get that $w$   has exactly   one  adjacent vertex  in $\{x_{4}, x_{5}\}$.
 Since  $u$, $v$ and  $w$  have  no  common adjacent vertex in  $N_{G}(x)$,  without loss of generality,  let  $v$  be adjacent  to  $x_{4}$  and  $w$  be adjacent  to  $x_{5}$.   But $G$  contains a $W_{4}$ with vertex set $\{u,v,w,x_4,x_5\}$,  contradicting \autoref{lemma:2.1}.
 %But  $G+vx_{5}$  has  two  $K_{4}$  with  vertex  sets  $\{u,v,x_{4},x_{5}\}$  and  $\{u,v,w,x_{5}\}$,   a contradiction.

If  $K_{4}$  in  $G+xu$ and  $K_{4}$  in  $G+xw$  contain  no vertex  in $\{x_{1}, x_{2}\}$,   then let  $V(K_{4})=\{x,w,x_{i},x_{j}\}$  in  $G+xw$  and  $V(K_{4})=\{x,u,x_{s},x_{t}\}$  in  $G+xu$,  and  $|\{x_{i},x_{j}\}\cap\{x_{s},x_{t}\}|\leq1$ since  $G$ is $K_4$-free.
If $|\{x_{i},x_{j}\}\cap\{x_{s},x_{t}\}|=1$,  then it is the same as the situation described in the paragraph 3 of Case 2 since $u,v,w$ are symmetric.  If $|\{x_{i},x_{j}\}\cap\{x_{s},x_{t}\}|=0$,  then, without loss of generality,   we say  $V(K_{4})=\{x,w,x_{3},x_{4}\}$  in  $G+xw$  and  $V(K_{4})=\{x,u,x_{5},x_{6}\}$  in  $G+xu$.   By a similar reasoning as that in the proof of paragraph 3 of Case 2,    $v$  and  $w$  have  exactly   one  adjacent vertex  in $\{x_{5}, x_{6}\}$.
Since  $u$, $v$ and  $w$  have  no  common adjacent vertex, without loss of generality,  we  say  that $v$  is  adjacent  to   $x_{5}$   and  $w$  is  adjacent  to   $x_{6}$.   But $G$  contains  a $W_{4}$ with vertex set $\{u,v,w,x_5,x_6\}$,  contradicting \autoref{lemma:2.1}.

Therefore,  by the above discussion,  there are no uniquely $K_{4}$-saturated graphs  on  $n$  vertices  with  maximum degree  $n-4$.  Our result follows immediately.
\end{proof}

\begin{theorem}\label{theorem:3.2}
If  $G$ is  a  uniquely $K_{4}$-saturated  graph  on  $n$  vertices  with  maximum degree  $n-3$,   then  $G=\overline{C_{7}}$.
\end{theorem}
\begin{proof}
Firstly we can get that $\overline{C_{7}}$  is a uniquely $K_{4}$-saturated  graph  by a routine check.

Let  $G$  be a  uniquely $K_{4}$-saturated  graph  on  $n$  vertices  with  maximum degree  $n-3$   and  $x$  be a vertex with degree  $n-3$.   Let  $N_{G}(x)=\{x_{1},x_{2},\ldots,x_{n-3}\}$   and $A=V(G)\setminus N_{G}[x]=\{y, z\}$.  $G[N_{G}(x)]$  is $K_{3}$-free  since  $G$  is  $K_{4}$-free.   By  \autoref{lemma:2.2},  we  can  get  that  $y$  is adjacent to $z$.
 $G+xy$  contains  exactly  one  $K_{4}$  and  $z\notin K_{4}$  since  $z$  is  not  adjacent to $x$.  Without loss of generality,  we  say   that  $K_{4}$ in $G+xy$  has  vertex set  $V(K_{4})=\{x,y,x_{1},x_{2}\}$.    Then $K_{4}$  in   $G+xz$  contains at most one vertex in $\{x_{1}, x_{2}\}$  since  $G$  is  $K_{4}$-free.

If  $K_{4}$  in  $G+xz$  contains exactly  one vertex in $\{x_{1}, x_{2}\}$,  then,  without loss of generality,  let $V(K_{4})=\{x,z,x_{2},x_{3}\}$  in  $G+xz$.
Now  we claim  $N_{G}(z)=\big(N_{G}(x)-\{x_{1}\}\big)\cup\{y\}$.  Otherwise,  for  $i\geq4$,  there  exists a vertex  $x_{i}$,  say  $x_{4}$,  which  is not adjacent to  $z$.   $G+zx_{4}$  contains  one $K_{4}$,   then  $y\in K_{4}$ and  $x\notin K_{4}$ since $G[N_{G}(x)]$ is  $K_{3}$-free  and  $x$  is not adjacent to  $z$.  Thus  $K_{4}$  in  $G+zx_{4}$  has  vertex set  $\{z,y,x_{4},x_{i}\}$  with  $i\neq1$.  But  $G+xy$  contains  two $K_{4}$  with  vertex sets  $\{x,y,x_{1},x_{2}\}$    and  $\{x,y,x_{4},x_{i}\}$,    a contradiction.
Therefore,  $N_{G}(z)=\big(N_{G}(x)-\{x_{1}\}\big)\cup\{y\}$.  Similarly,  we have $N_{G}(y)=\big(N_{G}(x)-\{x_{3}\}\big)\cup\{z\}$.
Moreover,  we claim  $G[N_{G}(x)]=P_{3}\cup(n-6)K_{1}$  with  $P_{3}=x_{1}x_{2}x_{3}$.   If not,   $G+xy$  or  $G+xz$  contains  more than one $K_{4}$,  a contradiction.  Moreover $|N_{G}(x)|\geq5$  since $\Delta=n-3$.  But  $G+x_{1}x_{4}$  contains  no  $K_{4}$, a contradiction  again.

If  $K_{4}$  in  $G+xz$  contains no vertex in $\{x_{1}, x_{2}\}$,  then,  without loss of generality,  let $V(K_{4})=\{x,z,x_{3},x_{4}\}$  in  $G+xz$.
Similarly as the proof of  previous paragraph, we can get that   $y$  and $z$  are  adjacent  to  $x_{i}$  for  any  $i\geq5$.
Moreover,  $y$ has  at most one  adjacent vertex  in $\{x_{3}, x_{4}\}$   and $z$ has  at most one  adjacent vertex  in $\{x_{1}, x_{2}\}$   since $G$  is  $K_{4}$-free.
If  $y$  is  not adjacent to $x_{3}$, $x_{4}$,   then  $G+yx_{3}$  contains one $K_{4}$.   We have   $z\in K_{4}$  and  $x, x_{4}\notin K_{4}$   in $G+yx_{3}$  since  $G[N_{G}(x)]$  is $K_{3}$-free  and  $x$, $x_{4}$  are  not adjacent to  $y$.   Then   $K_{4}$  in  $G+yx_{3}$  has  vertex set  $\{y,z,x_{3},x_{i}\}$  with  $i\neq4$.  But  $G+xz$  contains  two $K_{4}$  with  vertex sets  $\{z,x,x_{3},x_{4}\}$    and  $\{z,x,x_{3},x_{i}\}$,    a contradiction.  Thus  $y$  has  exactly  one  adjacent vertex  in $\{x_{3}, x_{4}\}$.   Similarly as above, we can also get that  $z$ has  exactly  one  adjacent vertex  in $\{x_{1}, x_{2}\}$.   Therefore,  without loss of generality,  we  say  $N_{G}(y)=\big(N_{G}(x)-\{x_{3}\}\big)\cup\{z\}$  and  $N_{G}(z)=\big(N_{G}(x)-\{x_{2}\}\big)\cup\{y\}$.
Since  $G+xy$  and  $G+xz$  contain exactly  one  $K_{4}$,  we can get  $E(N_{G}(x)\setminus\{x_{1},x_{2},x_{3},x_{4}\}, \{x_{1},x_{2},x_{3},x_{4}\})=\emptyset$,  $x_{1}$  is  not adjacent to  $x_{3}$ and  $x_{4}$,    $x_{2}$  is  not adjacent to  $x_{4}$.   Moreover,  we  claim  that  $x_{2}$  is adjacent to  $x_{3}$.  If  not,   $G+x_{2}x_{3}$  contains  no  $K_{4}$  by  the structure of $G$,   a contradiction.  Thus  $G[N_{G}(x)]=P_{4}\cup(n-7)K_{1}$  with  $P_{4}=x_{1}x_{2}x_{3}x_{4}$   and  $|N_{G}(x)|\geq4$.
If   $|N_{G}(x)|=4$,  we get $G=\overline{C_{7}}$ with $\overline{G}=C_{7}=xyx_3x_1x_4x_2zx$ by a routine check.   If   $|N_{G}(x)|\geq5$,   for any  $x_{i}\in N_{G}(x)\setminus\{x_{1},x_{2},x_{3},x_{4}\}$,  $G+x_{2}x_{i}$  contains  no $K_{4}$,  a contradiction.

Therefore,  $\overline{C_{7}}$  is the only  uniquely $K_{4}$-saturated  graph  on  $n$  vertices  with  maximum degree  $n-3$.
\end{proof}

\begin{center}
\begin{tikzpicture} [inner sep=0.8mm]
\tikzstyle{place}=[circle,draw=black!,fill=black!]
\tikzstyle{r}=[circle,draw=red!,fill=red!]

\node at (0,0)    [place] {};
\node at (0,1)    [place] {};
\node at (0,2)    [place] {};
\node at (0,3)    [place] {};
\node at (0,4)    [place] {};

\node at (5,0)    [place] {};
\node at (5,1)    [place] {};
\node at (5,2)    [place] {};
\node at (5,3)    [place] {};
\node at (5,4)    [place] {};

\draw[black,line width=1] (0,0) -- (0,1);
\draw[black,line width=1] (0,1) -- (0,2);
\draw[black,line width=1] (0,2) -- (0,3);
\draw[black,line width=1] (0,3) -- (0,4);
\draw[black,line width=1] (0,4) .. controls(-1,3)  and  (-1, 1)   .. (0,0);

\draw[black,line width=1] (5,0) -- (5,1);
\draw[black,line width=1] (5,1) -- (5,2);
\draw[black,line width=1] (5,2) -- (5,3);
\draw[black,line width=1] (5,3) -- (5,4);
\draw[black,line width=1] (5,4) .. controls(6,3)  and  (6, 1)   .. (5,0);

\draw[black,line width=1] (0,0) -- (5,0);
\draw[black,line width=1] (0,0) -- (5,1);
\draw[black,line width=1] (0,0) -- (5,2);

\draw[black,line width=1] (0,1) -- (5,2);
\draw[black,line width=1] (0,1) -- (5,3);
\draw[black,line width=1] (0,1) -- (5,4);

\draw[black,line width=1] (0,2) -- (5,0);
\draw[black,line width=1] (0,2) -- (5,1);
\draw[black,line width=1] (0,2) -- (5,4);

\draw[black,line width=1] (0,3) -- (5,1);
\draw[black,line width=1] (0,3) -- (5,2);
\draw[black,line width=1] (0,3) -- (5,3);

\draw[black,line width=1] (0,4) -- (5,0);
\draw[black,line width=1] (0,4) -- (5,3);
\draw[black,line width=1] (0,4) -- (5,4);

\coordinate [label=below left:$a_{4}$] () at (0,0);
\coordinate [label=above left:$a_{3}$] () at (0,1);
\coordinate [label= left:$a_{2}$] () at (0,2);
\coordinate [label=below left:$a_{1}$] () at (0,3);
\coordinate [label=above left:$a_{0}$] () at (0,4);

\coordinate [label=below right:$b_{4}$] () at (5,0);
\coordinate [label=above right:$b_{3}$] () at (5,1);
\coordinate [label= right:$b_{2}$] () at (5,2);
\coordinate [label=below right:$b_{1}$] () at (5,3);
\coordinate [label=above right:$b_{0}$] () at (5,4);

\end{tikzpicture}
\end{center}
\begin{center}
 Figure 2 :  $4$-primitive graph $G_{10}$.
\end{center}

In \autoref{theorem:3.3}, we prove   that  the  only   uniquely $K_{4}$-saturated  graph   on  $n$  vertices  with  maximum degree  $n-5$ is  $G_{10}$   in  Figure 2. The proof uses the similar approaches of  \autoref{theorem:3.1}  and  \autoref{theorem:3.2},  but  there  are   more classifications about  $E(G[V(G)\backslash N_{G}[x]])$  with  maximum  degree  vertex $x$  by  \autoref{lemma:2.2}.   More  details are needed to consider  to exclude  some  cases when $\Delta=n-5$.   We appropriately omit  some similar proofs  in  \autoref{theorem:3.3}.

\begin{theorem}\label{theorem:3.3}
 If  $G$ is  a  uniquely $K_{4}$-saturated  graph  on  $n$  vertices  with  maximum degree  $n-5$,   then  $G=G_{10}$   in  Figure 2.
\end{theorem}
\begin{proof}
Firstly we can get that $G_{10}$ in Figure 2 is a uniquely $K_{4}$-saturated  graph by a routine check.

Let  $G$  be a  uniquely $K_{4}$-saturated graph  on  $n$  vertices  with  maximum degree  $n-5$  and  $x$  be a vertex with degree  $n-5$.  Let   $N_{G}(x)=\{x_{1},x_{2},\ldots,x_{n-5}\}$   and    $A=V(G)\setminus N_{G}[x]=\{y_{1}, y_{2}, y_{3}, y_{4}\}$.  $G[N_{G}(x)]$  is  $K_{3}$-free   since  $G$  is  $K_{4}$-free.   By  \autoref{lemma:2.2}  and  $G$  is  $K_{4}$-free,  we have  $2\leq e(G[A])\leq5$.

\textbf{Case 1.} $e(G[A])=2$.

By  \autoref{lemma:2.2},   without loss of generality,  we  say that  $y_{1}$ is adjacent to $y_{2}$,  $y_{3}$ is adjacent to $y_{4}$.    $G+xy_{1}$  contains  exactly  one  $K_{4}$  and  $y_{2}, y_{3}, y_{4}\notin K_{4}$  since  $y_{2}, y_{3}, y_{4}$  are  not  adjacent to $x$.  Without loss of generality,  we  say   that  $K_{4}$ in $G+xy_{1}$  has  vertex set  $V(K_{4})=\{x,y_{1},x_{1},x_{2}\}$.    $K_{4}$  in   $G+xy_{2}$  contains at most one vertex in $\{x_{1}, x_{2}\}$  since  $G$  is  $K_{4}$-free.

If  $K_{4}$  in  $G+xy_{2}$  contains exactly  one vertex in $\{x_{1}, x_{2}\}$,  then, without loss of generality,  let  $V(K_{4})=\{x,y_{2},x_{2},x_{3}\}$  in  $G+xy_{2}$.
By a similar reasoning as that in the proof (paragraph 3)  of   \autoref{theorem:3.2},    we can get $N_{G}(y_{1})=\big(N_{G}(x)-\{x_{3}\}\big)\cup\{y_{2}\}$  and  $N_{G}(y_{2})=\big(N_{G}(x)-\{x_{1}\}\big)\cup\{y_{1}\}$.
Then  $G[N_{G}(x)]=P_{3}\cup(n-8)K_{1}$  with  $P_{3}=x_{1}x_{2}x_{3}$  since   $G+xy_{1}$  and  $G+xy_{2}$  contain  exactly  one $K_{4}$.    Thus by the  structure of  $G[N_{G}(x)]$,  we  have  that  $K_{4}$  in  $G+y_{1}y_{i}$  has vertex set $\{y_{1},y_{i},x_{1},x_{2}\}$   for  $i\in\{3,4\}$.   But  $G$  contains a $K_{4}$  with  vertex set  $\{y_{3},y_{4},x_{1},x_{2}\}$,  a contradiction.

If  $K_{4}$  in  $G+xy_{2}$  contains no vertex in $\{x_{1}, x_{2}\}$,  then, without loss of generality,  let   $V(K_{4})=\{x,y_{2},x_{3},x_{4}\}$  in  $G+xy_{2}$.   By a similar reasoning as that in the proof (paragraph 4)  of   \autoref{theorem:3.2}, without loss of generality,  we can get $N_{G}(y_{1})=\big(N_{G}(x)-\{x_{3}\}\big)\cup\{y_{2}\}$   and    $N_{G}(y_{2})=\big(N_{G}(x)-\{x_{2}\}\big)\cup\{y_{1}\}$.
 Then  $G[N_{G}(x)]=2P_{2}\cup(n-9)K_{1}$  with  $2P_{2}=x_{1}x_{2}\cup x_{3}x_{4}$   or  $G[N_{G}(x)]=P_{4}\cup(n-9)K_{1}$  with  $P_{4}=x_{1}x_{2}x_{3}x_{4}$  since   $G+xy_{1}$  and  $G+xy_{2}$  contain  exactly  one $K_{4}$.   Thus by the  structure of  $G[N_{G}(x)]$,  $K_{4}$  in  $G+y_{1}y_{i}$  has vertex set $\{y_{1},y_{i},x_{1},x_{2}\}$   for  $i\in\{3,4\}$.   But  $G$  contains a $K_{4}$  with  vertex set  $\{y_{3},y_{4},x_{1},x_{2}\}$,  a contradiction.

 Therefore,  $G$  is not  uniquely $K_{4}$-saturated  if  $e(G[A])=2$.

\textbf{Case 2.} $e(G[A])=3$.

By  \autoref{lemma:2.2},  we have $G[A]=P_4$  or $G[A]=K_{1,3}$.

\textbf{Subcase 2.1.} $G[A]=P_4$.

Without loss of generality,  we  say that  $y_{2}$ is adjacent to $y_{1}$  and  $y_{3}$,    $y_{3}$ is adjacent to $y_{4}$.  $G+xy_{2}$  contains exactly  one  $K_{4}$  and  $y_{1}, y_{3}, y_{4}\notin K_{4}$  since  $y_{1}, y_{3}, y_{4}$  are  not  adjacent to $x$.  Without loss of generality,  we  say   that  $K_{4}$ in $G+xy_{2}$  has  vertex set  $V(K_{4})=\{x,y_{2},x_{1},x_{2}\}$.  $K_{4}$  in   $G+xy_{3}$  contains at most one vertex in $\{x_{1}, x_{2}\}$  since  $G$  is  $K_{4}$-free.

If  $K_{4}$  in  $G+xy_{3}$  contains no vertex in $\{x_{1}, x_{2}\}$,  then, without loss of generality,  let   $V(K_{4})=\{x,y_{3},x_{3},x_{4}\}$  in  $G+xy_{3}$.   $K_{4}$  in   $G+xy_{1}$  contains at most one vertex in $\{x_{1}, x_{2}\}$   and   $K_{4}$  in   $G+xy_{4}$  contains at most one vertex in $\{x_{3}, x_{4}\}$  since $G$ is $K_4$-free.   By a similar reasoning as that in the proof (paragraph 3  and  4)  of   \autoref{theorem:3.2},  we can get
 $N_{G}(y_{1})=\big(N_{G}(x)-\{x_{i}\}\big)\cup\{y_{2}\}$  with  $i=1$  or  $2$  and    $N_{G}(y_{4})=\big(N_{G}(x)-\{x_{j}\}\big)\cup\{y_{3}\}$  with  $j=3$  or  $4$.
Without loss of generality,  we  say  $i=1$  and  $j=4$,  that is,  $N_{G}(y_{1})=\big(N_{G}(x)-\{x_{1}\}\big)\cup\{y_{2}\}$   and    $N_{G}(y_{4})=\big(N_{G}(x)-\{x_{4}\}\big)\cup\{y_{3}\}$.  Moreover,  $G+y_{1}y_{4}$  contains exactly  one  $K_{4}$  and  $y_{2}, y_{3}\notin K_{4}$,   then   $K_{4}$ in $G+y_{1}y_{4}$  has  vertex set  $V(K_{4})=\{y_{1},y_{4},x_{i},x_{j}\}$  with  $i,j\notin\{1,4\}$.   But   $G+xy_{1}$  contains  two $K_{4}$  with  vertex sets  $\{x,y_{1},x_{3},x_{4}\}$  and  $\{x,y_{1},x_{i},x_{j}\}$, a contradiction.

If  $K_{4}$  in  $G+xy_{3}$  contains exactly  one  vertex in $\{x_{1}, x_{2}\}$,  then, without loss of generality,  let  $V(K_{4})=\{x,y_{3},x_{2},x_{3}\}$  in  $G+xy_{3}$.   $K_{4}$  in   $G+xy_{1}$  contains at most one vertex in $\{x_{1}, x_{2}\}$   and   $K_{4}$  in   $G+xy_{4}$  contains at most one vertex in $\{x_{2}, x_{3}\}$ since $G$ is $K_4$-free.   By a similar reasoning as that in the proof (paragraph 3  and  4)  of   \autoref{theorem:3.2},  we can get
 $N_{G}(y_{1})=\big(N_{G}(x)-\{x_{i}\}\big)\cup\{y_{2}\}$  with  $i=1$  or  $2$  and    $N_{G}(y_{4})=\big(N_{G}(x)-\{x_{j}\}\big)\cup\{y_{3}\}$  with  $j=2$  or  $3$.
 If  $i=1$,  then  $N_{G}(y_{1})=\big(N_{G}(x)-\{x_{1}\}\big)\cup\{y_{2}\}$.  But $G$ contains a $W_4$ with vertex set $\{x_2,x_3,y_1,y_2,y_3\}$, contradicting \autoref{lemma:2.1}.
 %But  $G+y_{2}x_{3}$   has  two $K_{4}$  with  vertex sets  $\{y_{2},x_{3},y_{3},x_{2}\}$  and  $\{y_{2},x_{3},y_{1},x_{2}\}$, a contradiction.
 Thus $N_{G}(y_{1})=\big(N_{G}(x)-\{x_{2}\}\big)\cup\{y_{2}\}$.
 If  $j=3$,  then  $N_{G}(y_{4})=\big(N_{G}(x)-\{x_{3}\}\big)\cup\{y_{3}\}$.   But $G$ contains a $W_4$ with vertex set $\{x_1,x_2,y_2,y_3,y_4\}$, contradicting \autoref{lemma:2.1}.
 Thus $N_{G}(y_{4})=\big(N_{G}(x)-\{x_{2}\}\big)\cup\{y_{3}\}$.
 $G+y_{1}y_{4}$  contains  one $K_{4}$,  then  there exist  two adjacent  vertices  $x_{p}$  and  $x_{q}$ in $N_{G}(x)\setminus\{x_{2}\}$ such  that  $y_{i}$ is   adjacent to  $x_{p}$  and  $x_{q}$  with  $i=1, 4$.   Now  we claim  $p,q\notin\{1,3\}$.  Otherwise,  without loss of generality,  we  assume  $q=1$  and  $p=4$.  Moreover,  we have  that  $y_{3}$  is not adjacent to $x_{4}$.  If not,
 $G$ contains a $F^{-}(2,4,2)$ with vertex set $\{x_1,x_2,x_4,y_2,y_3,y_4\}$, contradicting \autoref{lemma:2.1}.
 Thus   $y_{3}$  is not adjacent to $x_{4}$.  Then
 $G+y_{3}x_{4}$  contains  one $K_{4}$  and  $y_{4}\in K_{4}$,  $y_{2}, x_{1}\notin K_{4}$   since  $G[N_{G}(x)]$  is  $K_{3}$-free  and  $G$  is  $K_{4}$-free.  Then   $K_{4}$  in  $G+y_{3}x_{4}$  has  vertex set  $\{y_{3},x_{4},y_{4},x_{i}\}$  with  $i\neq1$.
 But $G$ contains a $W_4$ with vertex set $\{x,x_1,x_4,x_i,y_4\}$, contradicting \autoref{lemma:2.1}.
 Thus   $p,q\notin\{1,3\}$,  without loss of generality,  we   assume    $p=4$  and  $q=5$.   Since  $G$  is $K_{4}$-free,  $y_{2}, y_{3}$  have at most one  adjacent vertex in  $\{x_{4}, x_{5}\}$.
 If  $y_{2}$  has no  adjacent vertex in  $\{x_{4}, x_{5}\}$,  then  $G+y_{2}x_{5}$  contains  one $K_{4}$  and  $x_{4}\notin K_{4}$.  Moreover,   $y_{1}\in K_{4}$  or  $y_{3}\in K_{4}$ in $G+y_{2}x_{5}$  since $G[N_G(x)]$ is $K_3$-free, but  this means that  $G+xy_{1}$   or   $G+xy_{3}$    contains  two $K_{4}$,  a contradiction.   Thus  $y_{2}$  has exactly  one  adjacent vertex in  $\{x_{4}, x_{5}\}$.  Similarly  as above,   we can also get that $y_{3}$  has exactly  one  adjacent vertex in  $\{x_{4}, x_{5}\}$.
  For  $i\geq6$,  if  there  exists $x_{i}$  such  that $x_{i}$  is  not  adjacent to  $y_{2}$,  then  $G+y_{2}x_{i}$  contains one  $K_{4}$   and  we  have  $y_{1}\in K_{4}$  or  $y_{3}\in K_{4}$    since  $G[N_{G}(x)]$  is  $K_{3}$-free.    But this means that  $G+xy_{1}$   or   $G+xy_{3}$    contains  two $K_{4}$,  a contradiction.  Thus  $y_{2}$  is adjacent  to  $x_{i}$  for  any  $i\geq6$. Similarly,  we have that $y_{3}$  is adjacent  to  $x_{i}$  for  any  $i\geq6$.  Moreover,  $y_{2}$  is  not adjacent  to $x_{3}$,  and  $y_{3}$  is  not adjacent  to $x_{1}$  since  $G$  is  $K_{4}$-free.   Therefore,  we  have
  $N_{G}(y_{2})=\big(N_{G}(x)-\{x_{3}, x_{i}\}\big)\cup\{y_{1}, y_{3}\}$  with  $i=4$  or $5$  and    $N_{G}(y_{3})=\big(N_{G}(x)-\{x_{1}, x_{j}\}\big)\cup\{y_{2}, y_{4}\}$  with  $j=4$  or $5$.
Moreover, we can claim $i\neq j$.  Otherwise, we assume $i=j=4$, then  $|N_{G}(x)|\geq6$  since  $\Delta\geq n-5$.  But $G+x_5x_6$ contains  more than one $K_{4}$,  a contradiction.  Thus $i\neq j$.
For $t\in[4]$,  since  $G+xy_{t}$    contains  exactly  one  $K_{4}$,   we  have  $G[N_{G}(x)]=P_{3}\cup P_{2}\cup(n-10)K_{1}$  with  $P_{3}=x_{1}x_{2}x_{3}$  and  $P_{2}=x_{4}x_{5}$.    But $G+x_{2}x_{4}$   or  $G+x_{2}x_{5}$   contains  no  $K_{4}$  by the  structure of $G$,  a contradiction.

\textbf{Subcase 2.2.} $G[A]=K_{1,3}$.

Without loss of generality,  we  say that  $y_{1}$ is adjacent to $y_{2}$,  $y_{3}$  and  $y_4$.
$G+xy_{1}$  contains exactly  one  $K_{4}$  and  $y_{2}, y_{3}, y_{4}\notin K_{4}$  since  $y_{2}, y_{3}, y_{4}$  are  not  adjacent to $x$.  Without loss of generality,  we  say   that  $K_{4}$ in $G+xy_{1}$  has  vertex set  $V(K_{4})=\{x,y_{1},x_{1},x_{2}\}$.  $K_{4}$  in   $G+xy_{i}$,  $i\in\{2,3,4\}$,  contains at most one vertex in $\{x_{1}, x_{2}\}$  since  $G$  is  $K_{4}$-free.
By a similar reasoning as that in the proof (paragraph 3  and  4)  of   \autoref{theorem:3.2},  we can get that for  $i\in\{2,3,4\}$, $N_{G}(y_{i})=\big(N_{G}(x)-\{x_{1}\}\big)\cup\{y_{1}\}$  or  $N_{G}(y_{i})=\big(N_{G}(x)-\{x_{2}\}\big)\cup\{y_{1}\}$.  Thus there are two vertices in $\{y_2,y_3,y_4\}$,  say $y_2,y_3$, such that $N_{G}(y_{2})=N_{G}(y_{3})$.  Without loss of generality,  let $N_{G}(y_{2})=N_{G}(y_{3})=\big(N_{G}(x)-\{x_{1}\}\big)\cup\{y_{1}\}$.  There exist  two adjacent  vertices  $x_{p}$  and  $x_{q}$ in $N_{G}(x)\setminus\{x_{1}\}$  since $G+xy_2$ contains a $K_4$.  But  $G+y_2y_3$ contains two $K_4$  with  vertex sets  $\{x_2,y_{1},y_{2},y_{3}\}$  and  $\{y_2,y_{3},x_{p},x_{q}\}$, a contradiction.

 Therefore,  $G$  is  not uniquely $K_4$-saturated  if  $e(G[A])=3$.

\textbf{Case 3.} $e(G[A])=4$.

\textbf{Subcase 3.1.} $G[A]$  is formed by adding a  pendant edge to one vertex on a triangle.

Without loss of generality,  we  say that  $y_{1}$ is adjacent to $y_{2}$  and  $y_{3}$,  $y_{2}$ is adjacent to $y_{3}$,  $y_{3}$ is adjacent to $y_{4}$.  Moreover,  $y_{1}$, $y_{2}$ and  $y_{3}$  have  no  common adjacent vertex in  $N_{G}(x)$  since  $G$ is  $K_{4}$-free.
$G+xy_{2}$  contains  exactly  one  $K_{4}$  and  $y_{1}, y_{3}, y_{4}\notin K_{4}$  since  $y_{1}, y_{3}, y_{4}$  are  not  adjacent to $x$.  Without loss of generality,  we  say   that  $K_{4}$ in $G+xy_{2}$  has  vertex set  $V(K_{4})=\{x,y_{2},x_{1},x_{2}\}$.  Then $K_{4}$  in   $G+xy_{1}$  and  $K_4$  in   $G+xy_{3}$  contain at most one vertex in $\{x_{1}, x_{2}\}$  since  $G$  is  $K_{4}$-free.

If  $K_{4}$  in  $G+xy_{1}$  and  $K_4$ in  $G+xy_{3}$   contain exactly  one vertex in $\{x_{1}, x_{2}\}$,  then,  without loss of generality,    let  $V(K_{4})=\{x,y_{1},x_{1},x_{4}\}$  in  $G+xy_{1}$  and  $V(K_{4})=\{x,y_{3},x_{2},x_{3}\}$  in  $G+xy_{3}$.   $y_{3}$  is  not  adjacent to  $x_{1}$  since $G$  is  $K_{4}$-free.
But $G$ contains a $W_4$ with vertex set $\{x_1,x_2,y_1,y_2,y_3\}$, contradicting \autoref{lemma:2.1}.

If $K_{4}$  in   $G+xy_{3}$  contains exactly  one vertex in $\{x_{1}, x_{2}\}$  and  $K_{4}$  in   $G+xy_{1}$  contains no vertex in $\{x_{1}, x_{2}\}$, then, without loss of generality, let $V(K_{4})=\{x,y_{3},x_{2},x_{3}\}$  in  $G+xy_{3}$.
If  $K_{4}$  in  $G+xy_{1}$  and  $K_4$  in $G+xy_{3}$ have a common vertex  in  $N_{G}(x)$,  then we assume that  $V(K_{4})=\{x,y_{1},x_{3},x_{4}\}$  in  $G+xy_{1}$. 
But $G$ contains a $W_4$ with vertex set $\{x_2,x_3,y_1,y_2,y_3\}$, contradicting \autoref{lemma:2.1}.
If  $K_{4}$  in  $G+xy_{1}$  and  $K_4$ in $G+xy_{3}$ have  no  common vertex  in  $N_{G}(x)$,  then, without loss of generality,   let  $V(K_{4})=\{x,y_{3},x_{2},x_{3}\}$  in  $G+xy_{3}$   and   $V(K_{4})=\{x,y_{1},x_{4},x_{5}\}$  in  $G+xy_{1}$.  Moreover,  $y_{1}$  is  not  adjacent  to  $x_{2}$,   $y_{3}$  is  not  adjacent  to  $x_{1}$  since $G$  is  $K_{4}$-free. We claim that  $y_{1}$  is  not  adjacent  to  $x_{1}$,  if  not,  $G$ contains a $W_4$ with vertex set $\{x_1,x_2,y_1,y_2,y_3\}$, contradicting \autoref{lemma:2.1}.
Then  $G+y_{1}x_{1}$   contains  one  $K_{4}$  and
$x_{2}, y_{3}\notin K_{4}$,  $y_{2}\in K_{4}$   since  $G[N_{G}(x)]$  is  $K_{3}$-free.  But this means that $G+xy_{2}$    contains  two $K_{4}$,  a contradiction.

If  $K_{4}$  in   $G+xy_{1}$  contains exactly  one vertex in $\{x_{1}, x_{2}\}$  and  $K_{4}$  in   $G+xy_{3}$  contains no vertex in $\{x_{1}, x_{2}\}$,
then, without loss of generality, let  $V(K_{4})=\{x,y_{1},x_{1},x_{3}\}$  in  $G+xy_{1}$.
If  $K_{4}$  in  $G+xy_{1}$  and  $K_4$  in  $G+xy_{3}$  have a common vertex  in  $N_{G}(x)$,
then   we assume that  $V(K_{4})=\{x,y_{3},x_{3},x_{4}\}$  in  $G+xy_{3}$. 
But $G$ contains a $W_4$ with vertex set $\{x_1,x_3,y_1,y_2,y_3\}$, contradicting \autoref{lemma:2.1}.
If  $K_{4}$  in  $G+xy_{1}$  and $K_4$  in $G+xy_{3}$ have  no  common vertex  in  $N_{G}(x)$,   then, without loss of generality,   let  $V(K_{4})=\{x,y_{1},x_{1},x_{3}\}$  in  $G+xy_{1}$   and   $V(K_{4})=\{x,y_{3},x_{4},x_{5}\}$  in  $G+xy_{3}$.
Since  $G$  is  $K_{4}$-free,  $y_{1}$  and  $y_{2}$  have at most one  adjacent vertex  in $\{x_{4}, x_{5}\}$.
Now  we claim  that  $y_{1}$, $y_{2}$   have  exactly  one  adjacent vertex in  $\{x_{4}, x_{5}\}$.
Indeed, if $y_{2}$  is not adjacent to  $x_{4}$  and  $x_{5}$,  then  $G+y_{2}x_{4}$  contains  one $K_{4}$   and  $x_{5}\notin K_{4}$,  $y_{1}\in K_{4}$    since  $G[N_{G}(x)]$  is  $K_{3}$-free  and  $G+y_{3}x$  contains  exactly  one  $K_{4}$.   Thus $y_{1}$  is adjacent to  $x_{4}$.  Similarly,   we can  get  that  $G+y_{2}x_{5}$  contains  one $K_{4}$  and   $y_{1}\in K_{4}$,  which  means that  $y_{1}$  is adjacent to  $x_{5}$.  But  $G$  contains  a  $K_{4}$  with vertex  set  $\{y_{1},y_{3},x_{4},x_{5}\}$,  a contradiction.  Therefore,  $y_{2}$    has exactly   one  adjacent vertex  in $\{x_{4}, x_{5}\}$.   Similarly  as above proof, we can also get that $y_{1}$   has exactly   one  adjacent vertex  in $\{x_{4}, x_{5}\}$.   Since  $y_{1}$, $y_{2}$ and  $y_{3}$ have  no  common adjacent vertex in  $N_{G}(x)$,   without loss of generality,  we  assume  that  $y_{1}$  is adjacent  to  $x_{5}$,  $y_{2}$  is adjacent  to  $x_{4}$.   But $G$ contains a $W_4$ with vertex set $\{x_4,x_5,y_1,y_2,y_3\}$, contradicting \autoref{lemma:2.1}.
%But  $G+y_{2}x_{5}$  has  two  $K_{4}$  with  vertex  sets  $\{y_{2},y_{3},x_{4},x_{5}\}$  and  $\{y_{1},y_{2},y_{3},x_{5}\}$,   a contradiction.

If  $K_{4}$  in  $G+xy_{1}$ and $K_4$  in  $G+xy_{3}$  contain  no vertex in $\{x_{1}, x_{2}\}$,
 then, without loss of generality,   let $V(K_{4})=\{x,y_{3},x_{3},x_{4}\}$  in  $G+xy_{3}$.
If  $K_{4}$  in  $G+xy_{1}$  and  $K_4$  in  $G+xy_{3}$ have a common vertex  in  $N_{G}(x)$,
then we assume that $V(K_{4})=\{x,y_{1},x_{3},x_{5}\}$  in  $G+xy_{1}$.
Moreover,  $y_{2}$  is  not  adjacent  to  $x_{3}$,   $y_{3}$  is  not  adjacent  to  $x_{5}$   since $G$  is  $K_{4}$-free.  We claim that  $y_{2}$  is  not  adjacent  to  $x_{5}$,  if  not,  $G$ contains a $W_4$ with vertex set $\{x_3,x_5,y_1,y_2,y_3\}$, contradicting \autoref{lemma:2.1}.
Then  $G+y_{2}x_{5}$   contains  one  $K_{4}$  and
$x_{3}, y_{3}\notin K_{4}$,  $y_{1}\in K_{4}$  since  $G[N_{G}(x)]$  is  $K_{3}$-free.   But  $G+xy_{1}$    contains  two $K_{4}$,  a contradiction.
If  $K_{4}$  in  $G+xy_{1}$  and  $K_4$  in $G+xy_{3}$ have  no  common vertex  in  $N_{G}(x)$,  then, without loss of generality,   let $V(K_{4})=\{x,y_{3},x_{3},x_{4}\}$  in  $G+xy_{3}$   and   $V(K_{4})=\{x,y_{1},x_{5},x_{6}\}$  in  $G+xy_{1}$.
Since  $G$  is  $K_{4}$-free,  $y_{1}$  and  $y_{2}$  have at most one  adjacent vertex  in $\{x_{3}, x_{4}\}$.  Similarly  as the proof in the previous paragraph,  we get that  $y_{1}$ and $y_{2}$   have  exactly  one  adjacent vertex in  $\{x_{3}, x_{4}\}$.
Since  $y_{1}$, $y_{2}$ and  $y_{3}$ have  no  common adjacent vertex in  $N_{G}(x)$,   without loss of generality,  we  assume  that  $y_{1}$  is adjacent  to  $x_{3}$,  $y_{2}$  is adjacent  to  $x_{4}$.   But $G$ contains a $W_4$ with vertex set $\{x_3,x_4,y_1,y_2,y_3\}$, contradicting \autoref{lemma:2.1}.

\textbf{Subcase 3.2.} $G[A]=C_{4}$.

Without loss of generality,  we  say that  $y_{1}$ is adjacent to $y_{2}$  and  $y_{4}$,  $y_{2}$ is adjacent to $y_{3}$,  $y_{3}$ is adjacent to $y_{4}$.  Since  $G$  is  $K_{4}$-free,    any two adjacent  vertices  $y_{i}$,  $y_{j}$  have no  common adjacent edge  in  $G[N_{G}(x)]$.   We  divide into 35  cases  based  on the   number   of  vertices    in  $N_{G}(x)$,  say  $\alpha$,  which belong to $K_{4}$  in  $G+xy_{i}$,  $i\in[4]$. Then  $3\leq\alpha\leq8$  since  $G$  is  $K_{4}$-free.  We  get only  one  uniquely $K_4$-saturated graph in  this  subcase.
For convenience,  we  write  $y_{i}\rightarrow x_{p}, x_{q}$ to mean that the $K_{4}$  in  $G+xy_{i}$  has  vertex set $\{x, y_{i}, x_{p}, x_{q}\}$.

If  $\alpha=3$,  by the symmetry  of  $y_{i}$,  $i\in[4]$,   without loss of generality,  we can get that there are only  two  cases:   $y_{2}\rightarrow x_{1}, x_{2}$;  $y_{3}\rightarrow x_{2}, x_{3}$;  $y_{1}\rightarrow x_{2}, x_{3}$;   $y_{4}\rightarrow x_{1}, x_{2}$   and   $y_{2}\rightarrow x_{1}, x_{2}$;  $y_{3}\rightarrow x_{2}, x_{3}$;  $y_{1}\rightarrow x_{1}, x_{3}$;    $y_{4}\rightarrow x_{1}, x_{2}$.
But the  first case means  that $G$ contains a $W_4$ with vertex set $\{x_2,x_3,y_1,y_2,y_3\}$, contradicting \autoref{lemma:2.1}.
The  second   case  turns out  that  $G$  contains a $K_{4}$  with  vertex set  $\{x,x_{1},x_{2},x_{3}\}$,  a contradiction  again.

If  $\alpha=4$,  then  there are  13  cases  by the symmetry  of  $y_{i}$,  $i\in[4]$.   Without loss of generality,  we  say  that they are:
(a).  $y_{2}\rightarrow x_{1}, x_{2}$;  $y_{3}\rightarrow x_{2}, x_{3}$; $y_{1}\rightarrow x_{1}, x_{4}$;  $y_{4}\rightarrow x_{3}, x_{1}$
(b). $y_{2}\rightarrow x_{1}, x_{2}$;  $y_{3}\rightarrow x_{2}, x_{3}$; $y_{1}\rightarrow x_{1}, x_{4}$;    $y_{4}\rightarrow x_{3}, x_{4}$
(c). $y_{2}\rightarrow x_{1}, x_{2}$;  $y_{3}\rightarrow x_{2}, x_{3}$; $y_{1}\rightarrow x_{1}, x_{4}$;    $y_{4}\rightarrow x_{2}, x_{1}$
(d). $y_{2}\rightarrow x_{1}, x_{2}$;  $y_{3}\rightarrow x_{2}, x_{3}$; $y_{1}\rightarrow x_{1}, x_{4}$;    $y_{4}\rightarrow x_{2}, x_{4}$
(e). $y_{2}\rightarrow x_{1}, x_{2}$;  $y_{3}\rightarrow x_{2}, x_{3}$; $y_{1}\rightarrow x_{2}, x_{4}$;    $y_{4}\rightarrow x_{3}, x_{1}$
(f).  $y_{2}\rightarrow x_{1}, x_{2}$;  $y_{3}\rightarrow x_{2}, x_{3}$; $y_{1}\rightarrow x_{2}, x_{4}$;    $y_{4}\rightarrow x_{3}, x_{4}$
(g).  $y_{2}\rightarrow x_{1}, x_{2}$;  $y_{3}\rightarrow x_{2}, x_{3}$; $y_{1}\rightarrow x_{2}, x_{4}$;    $y_{4}\rightarrow x_{2}, x_{1}$
(h).  $y_{2}\rightarrow x_{1}, x_{2}$;  $y_{3}\rightarrow x_{2}, x_{3}$; $y_{1}\rightarrow x_{2}, x_{4}$;    $y_{4}\rightarrow x_{1}, x_{4}$
(i).  $y_{2}\rightarrow x_{1}, x_{2}$;  $y_{3}\rightarrow x_{2}, x_{3}$; $y_{1}\rightarrow x_{3}, x_{4}$;    $y_{4}\rightarrow x_{3}, x_{1}$
(j).  $y_{2}\rightarrow x_{1}, x_{2}$;  $y_{3}\rightarrow x_{2}, x_{3}$; $y_{1}\rightarrow x_{3}, x_{4}$;    $y_{4}\rightarrow x_{2}, x_{1}$
(k).  $y_{2}\rightarrow x_{1}, x_{2}$;  $y_{3}\rightarrow x_{2}, x_{3}$;  $y_{1}\rightarrow x_{3}, x_{4}$;    $y_{4}\rightarrow x_{2}, x_{4}$
(l).  $y_{2}\rightarrow x_{1}, x_{2}$;  $y_{3}\rightarrow x_{2}, x_{3}$;  $y_{1}\rightarrow x_{3}, x_{4}$;    $y_{4}\rightarrow x_{1}, x_{4}$
(m).  $y_{2}\rightarrow x_{1}, x_{2}$;  $y_{3}\rightarrow x_{3}, x_{4}$;  $y_{1}\rightarrow x_{3}, x_{4}$;    $y_{4}\rightarrow x_{1}, x_{2}$.
Cases  (a), (d), (e), (f), (h), (i)  and  (k)   mean  that  $G$  contains  a $K_{4}$,  contradicting the fact that  $G$  is  $K_{4}$-free.
For  cases  (b)  and  (l),   $G$ contains a $W_4$ with vertex set $\{x,x_1,x_2,x_3,x_4\}$, contradicting \autoref{lemma:2.1}.
For  cases  (c)  and  (j),   $G$ contains a $W_4$ with vertex set $\{x_1,x_2,y_2,y_3,y_4\}$, contradicting \autoref{lemma:2.1}.
For case (g),    $G$ contains a $W_4$ with vertex set $\{x_2,y_1,y_2,y_3,y_4\}$, contradicting \autoref{lemma:2.1}.
For case (m),  $y_{1}$  has   at most one adjacent vertex in $\{x_{1}, x_{2}\}$  since  $G$  is  $K_{4}$-free.   If  $y_{1}$  is  not adjacent to  $x_{1}$  and  $x_{2}$,  then  $G+y_{1}x_{1}$  contains  one  $K_{4}$  and  $x_{2}\notin K_{4}$.   Moreover,  we have   $y_{2}\in K_{4}$  or  $y_{4}\in K_{4}$  in  $G+y_{1}x_{1}$  since  $G[N_{G}(x)]$  is $K_3$-free.  But  $G+xy_{2}$  or  $G+xy_{4}$  contains  two  $K_{4}$  since  $x_{2}\notin K_{4}$,   a contradiction.  Thus   $y_{1}$  has   exactly  one adjacent vertex in $\{x_{1}, x_{2}\}$.   Without loss of generality,  we  assume  that  $y_{1}$ is adjacent to $x_{1}$.  But  $G$ contains a $W_4$ with vertex set $\{x_1,x_2,y_1,y_2,y_4\}$, contradicting \autoref{lemma:2.1}.

If  $\alpha=5$,  then  there are  11  cases  by the symmetry  of  $y_{i}$,  $i\in[4]$.   Without loss of generality,  we  say  that they are:
(a). $y_{2}\rightarrow x_{1}, x_{2}$;  $y_{3}\rightarrow x_{2}, x_{3}$; $y_{4}\rightarrow x_{3}, x_{5}$;    $y_{1}\rightarrow x_{1}, x_{4}$
(b).  $y_{2}\rightarrow x_{1}, x_{2}$;  $y_{3}\rightarrow x_{2}, x_{3}$; $y_{4}\rightarrow x_{3}, x_{5}$;    $y_{1}\rightarrow x_{2}, x_{4}$
(c). $y_{2}\rightarrow x_{1}, x_{2}$;  $y_{3}\rightarrow x_{2}, x_{3}$; $y_{4}\rightarrow x_{3}, x_{5}$;    $y_{1}\rightarrow x_{3}, x_{4}$
(d). $y_{2}\rightarrow x_{1}, x_{2}$;  $y_{3}\rightarrow x_{2}, x_{3}$;  $y_{4}\rightarrow x_{2}, x_{5}$;    $y_{1}\rightarrow x_{2}, x_{4}$
(e). $y_{2}\rightarrow x_{1}, x_{2}$;  $y_{3}\rightarrow x_{2}, x_{3}$; $y_{4}\rightarrow x_{2}, x_{5}$;    $y_{1}\rightarrow x_{3}, x_{4}$
(f).  $y_{2}\rightarrow x_{1}, x_{2}$;  $y_{3}\rightarrow x_{2}, x_{3}$; $y_{4}\rightarrow x_{2}, x_{5}$;    $y_{1}\rightarrow x_{5}, x_{4}$
(g). $y_{2}\rightarrow x_{1}, x_{2}$;  $y_{3}\rightarrow x_{2}, x_{3}$; $y_{4}\rightarrow x_{1}, x_{5}$;    $y_{1}\rightarrow x_{4}, x_{3}$
(h).  $y_{2}\rightarrow x_{1}, x_{2}$;  $y_{3}\rightarrow x_{2}, x_{3}$; $y_{4}\rightarrow x_{1}, x_{5}$;    $y_{1}\rightarrow x_{4}, x_{5}$
(i). $y_{2}\rightarrow x_{1}, x_{2}$;  $y_{3}\rightarrow x_{2}, x_{3}$; $y_{4}\rightarrow x_{1}, x_{3}$;    $y_{1}\rightarrow x_{4}, x_{5}$
(j). $y_{2}\rightarrow x_{1}, x_{2}$;  $y_{3}\rightarrow x_{2}, x_{3}$; $y_{4}\rightarrow x_{1}, x_{2}$;    $y_{1}\rightarrow x_{4}, x_{5}$
(k).  $y_{2}\rightarrow x_{1}, x_{2}$;  $y_{3}\rightarrow x_{3}, x_{4}$; $y_{4}\rightarrow x_{2}, x_{5}$;    $y_{1}\rightarrow x_{3}, x_{4}$.
For case (a),  $y_{2}$  is  not  adjacent to $x_{3}$  and $x_{4}$  since  $G$  is  $K_{4}$-free.  If  $y_{2}$  is  not adjacent to $x_{5}$,  then  $G+y_{2}x_{5}$  contains  one  $K_{4}$.  Moreover,   we have
$V(K_4)=\{y_1,y_2,x_5,x_i\}$  or  $V(K_4)=\{y_2,y_3,x_5,x_i\}$ with $i\neq3, 4$ in $G+y_{2}x_{5}$  since   $G[N_{G}(x)]$  is  $K_{3}$-free.   But then $G+xy_{1}$   or  $G+xy_{3}$   contains  two  $K_{4}$,   a contradiction.   Thus  $y_{2}$  is  adjacent to $x_{5}$.  Similarly  as  above,  we can get that $y_{3}$  is  adjacent to $x_{4}$,   $y_{3}$  is  not  adjacent to $x_{1}$ and $x_{5}$.
Since  $G$  is  $K_{4}$-free, $y_1$ is not adjacent to $x_2$,  and  $y_{1}$  has   at most one adjacent vertex in $\{x_{3}, x_{5}\}$.  If  $y_{1}$  is  not adjacent  to  $x_{3}$  and  $x_{5}$,  then   $G+y_{1}x_{3}$  contains  one  $K_{4}$  and  $x_{5}\notin K_{4}$.  Moreover, we have
$V(K_4)=\{y_1,y_2,x_3,x_i\}$  or  $V(K_4)=\{y_1,y_4,x_3,x_i\}$ with $i\neq5$ in $G+y_{1}x_{3}$  since   $G[N_{G}(x)]$  is  $K_{3}$-free.  But  $G+xy_{2}$   or  $G+xy_{4}$   contains  two  $K_{4}$,   a contradiction.  Thus  $y_{1}$  has  exactly  one adjacent vertex in $\{x_{3}, x_{5}\}$.
Similarly, we can get that  $y_{4}$ is not adjacent to $x_2$,  and  $y_{4}$  has  exactly  one adjacent vertex in $\{x_{1}, x_{4}\}$.
If  $y_{1}$  is adjacent to $x_{3}$  and  $y_{4}$  is adjacent to  $x_{1}$,  then  $G$ contains a $F^{-}(2,4,2)$ with vertex set $\{x,x_1,x_2,x_3,y_1,y_4\}$, contradicting \autoref{lemma:2.1}.
If  $y_{1}$  is adjacent to $x_{3}$  and  $y_{4}$  is adjacent to $x_{4}$,  then  $G$ contains a $W_4$ with vertex set $\{x_3,x_4,y_1,y_3,y_4\}$, contradicting \autoref{lemma:2.1}.
If  $y_{1}$  is adjacent to $x_{5}$  and  $y_{4}$  is adjacent to $x_{1}$,  then  $G$ contains a $W_4$ with vertex set $\{x_1,x_5,y_1,y_2,y_4\}$, contradicting \autoref{lemma:2.1}.
Thus  $y_{1}$  is adjacent to $x_{5}$  and  $y_{4}$  is adjacent to  $x_{4}$.
Moreover, we  have $G[\{x_{i}: i\in[5]\}]=P_{5}=x_4x_1x_2x_3x_5$  since  $G$ is $K_4$-free and  $G+xy_i$, $i\in[4]$,  contains exactly one $K_4$.
If  $|N_{G}(x)|=5$,  which  means $N_{G}(x)=\{x_{i}: i\in[5]\}$, then we get $G=G_{10}$  in Figure 2 with $a_0=x, a_1=x_4,  a_2=y_3, a_3=y_2,  a_4=x_5,  b_0=x_2,  b_1=x_1,  b_2=y_1,  b_3=y_4,  b_4=x_3$  by a routine check.
If   $|N_{G}(x)|\geq6$,  we  claim  that for  any  $i\geq6$  and  $j\in[4]$,  $y_{j}$  is  adjacent  to  $x_{i}$.   Otherwise,  without loss of generality,  we  say  that  $y_{1}$  is not  adjacent  to  $x_{i}$.  Then   $G+y_{1}x_{i}$  contains   one  $K_{4}$,  and
$y_{2}\in K_{4}$  or  $y_{4}\in K_{4}$  since   $G[N_{G}(x)]$  is  $K_{3}$-free.   But  $G+xy_{2}$   or  $G+xy_{4}$   contains  two  $K_{4}$,   a contradiction.  Thus for  $i\geq6$  and  $j\in[4]$,  $y_{j}$  is  adjacent  to  $x_{i}$.  And  $|N_{G}(x)|=6$ since $G+y_1y_3$  contains one $K_4$.
Moreover, $G[N_{G}(x)]=P_{5}\cup K_{1}$  since  $G+xy_{i}$, $i\in[4]$,  contains  exactly  one  $K_{4}$.  But  $G$ contains a $W_4$ with vertex set $\{x_4,x_6,y_1,y_3,y_4\}$, contradicting \autoref{lemma:2.1}.

For cases (b) and (d),  if  $y_{3}$  is  not adjacent to $x_{4}$,  then  $G+y_{3}x_{4}$  contains  one  $K_{4}$.  Moreover,   we have
$V(K_4)=\{y_3,y_2,x_4,x_i\}$  or  $V(K_4)=\{y_3,y_4,x_4,x_i\}$  in $G+y_{3}x_{4}$  since   $G[N_{G}(x)]$  is  $K_{3}$-free.   But then $G+xy_{2}$   or  $G+xy_{4}$   contains  two  $K_{4}$,   a contradiction.   Thus  $y_{3}$  is  adjacent to $x_{4}$.   But  $G$ contains a $W_4$ with vertex set $\{x,x_2,x_3,x_4,y_3\}$, contradicting \autoref{lemma:2.1}.
For cases (c) and (e),   if  $y_{3}$  is  not adjacent to $x_{4}$,  then $G+y_{3}x_{4}$  contains  one  $K_{4}$.  Moreover,   we have
$V(K_4)=\{y_3,y_2,x_4,x_i\}$  or  $V(K_4)=\{y_3,y_4,x_4,x_i\}$  in $G+y_{3}x_{4}$  since   $G[N_{G}(x)]$  is  $K_{3}$-free.   But then $G+xy_{2}$   or  $G+xy_{4}$   contains  two  $K_{4}$,   a contradiction.   Thus  $y_{3}$  is  adjacent to $x_{4}$.   But  $G$ contains a $W_4$ with vertex set $\{x,x_2,x_3,x_4,y_3\}$, contradicting \autoref{lemma:2.1}.
For case (f),  if  $y_{4}$  is  not adjacent to $x_{1}$,  then  $G+y_{4}x_{1}$  contains  one  $K_{4}$.  Moreover,   we have
$V(K_4)=\{y_4,y_3,x_1,x_i\}$  or  $V(K_4)=\{y_4,y_1,x_1,x_i\}$  in $G+y_{4}x_{1}$  since   $G[N_{G}(x)]$  is  $K_{3}$-free.   But then $G+xy_{3}$   or  $G+xy_{1}$   contains  two  $K_{4}$,   a contradiction.   Thus  $y_{4}$  is  adjacent to $x_{1}$.   But  $G$ contains a $W_4$ with vertex set $\{x_1,x_2,y_2,y_3,y_4\}$, contradicting \autoref{lemma:2.1}.
For case (g),  $y_1$  has  at most one adjacent vertex in  $\{x_1,x_2\}$  since  $G$  is $K_4$-free.   If  $y_{1}$  is  not adjacent  to  $x_{1}$  and  $x_{2}$,  then   $G+y_{1}x_{2}$  contains  one  $K_{4}$  and  $x_{1}\notin K_{4}$.  Moreover, we have
$V(K_4)=\{y_1,y_2,x_2,x_i\}$  or  $V(K_4)=\{y_1,y_4,x_2,x_i\}$ with $i\neq1$ in $G+y_{1}x_{2}$  since   $G[N_{G}(x)]$  is $K_{3}$-free.  But  $G+xy_{2}$   or  $G+xy_{4}$   contains  two  $K_{4}$,   a contradiction.  Thus  $y_{1}$  has  exactly  one adjacent vertex in $\{x_{1}, x_{2}\}$.  Similarly, we can get that  $y_{4}$  has  exactly  one adjacent vertex in $\{x_{2}, x_{3}\}$.
If  $y_{1}$  is adjacent to $x_{1}$  and  $y_{4}$  is adjacent to $x_{2}$,  then  $G$ contains a $F^{-}(2,4,2)$ with vertex set $\{x_1,x_2,y_1,y_2,y_3,y_4\}$, contradicting \autoref{lemma:2.1}.
If  $y_{1}$  is adjacent to $x_{1}$  and  $y_{4}$  is adjacent to $x_{3}$,  then  $G$ contains a $F^{-}(2,4,2)$ with vertex set $\{x,x_1,x_2,x_3,y_1,y_4\}$, contradicting \autoref{lemma:2.1}.
If  $y_{1}$  is adjacent to  $x_{2}$  and  $y_{4}$  is adjacent to $x_{2}$,  then  $G$ contains a $W_4$ with vertex set $\{x_2,y_1,y_2,y_3,y_4\}$, contradicting \autoref{lemma:2.1}.
If  $y_{1}$  is adjacent to $x_{2}$  and  $y_{4}$  is adjacent to $x_{3}$,  then  $G$ contains a $F^{-}(2,4,2)$ with vertex set $\{x_2,x_3,y_1,y_2,y_3,y_4\}$, contradicting \autoref{lemma:2.1}.
Thus  (g) does not hold.
For case (h), similar to  case (b),  we can get that $y_{1}$  is  adjacent to $x_{3}$.  $y_1$  is not adjacent to $x_{1}$  since $G$  is $K_4$-free.  Now we claim that $y_1$  is adjacent to $x_{2}$.  Otherwise, $G+y_1x_2$ contains one $K_4$ and  $x_1\notin K_4$.  Moreover,  we have
$V(K_4)=\{y_1,y_2,x_2,x_i\}$  or  $V(K_4)=\{y_1,y_4,x_2,x_i\}$  with  $i\neq1$  in $G+y_{1}x_{2}$  since   $G[N_{G}(x)]$  is  $K_{3}$-free.   But then $G+xy_{2}$   or  $G+xy_{4}$   contains  two  $K_{4}$,   a contradiction.   Thus  $y_{1}$  is  adjacent to $x_{2}$.   But $G$  contains a $W_4$ with vertex set $\{x_2,x_3,y_1,y_2,y_3\}$, contradicting \autoref{lemma:2.1}.
Case  (i)  means  that  $G$  contains  a $K_{4}$,  contradicting the fact that  $G$  is  $K_{4}$-free.
For  case (j),  $G$ contains a $W_4$ with vertex set $\{x_1,x_2,y_2,y_3,y_4\}$, contradicting \autoref{lemma:2.1}.
For  case (k),    $y_{2}$  has   at most one adjacent vertex in $\{x_{3}, x_{4}\}$  since  $G$  is  $K_{4}$-free.
If  $y_{2}$  is  not adjacent  to  $x_{3}$  and  $x_{4}$,  then   $G+y_{2}x_{3}$  contains  one  $K_{4}$  and  $x_{4}\notin K_{4}$.  Moreover, we have
$V(K_4)=\{y_1,y_2,x_3,x_i\}$  or  $V(K_4)=\{y_2,y_3,x_3,x_i\}$ with $i\neq4$ in $G+y_{2}x_{3}$  since   $G[N_{G}(x)]$  is $K_{3}$-free.  But  $G+xy_{1}$   or  $G+xy_{3}$   contains  two  $K_{4}$,   a contradiction.  Thus  $y_{2}$  has  exactly  one adjacent vertex in $\{x_{3}, x_{4}\}$.  Without loss of generality,  we assume that $y_{2}$  is adjacent to $x_{3}$.  But  $G$ contains a $W_4$ with vertex set $\{x_3,x_4,y_1,y_2,y_3\}$, contradicting \autoref{lemma:2.1}.

If  $\alpha=6$,  then  there are  6  cases  by the symmetry  of  $y_{i}$,  $i\in[4]$.   Without loss of generality,  we  say  that they are:
(a). $y_{2}\rightarrow x_{1}, x_{2}$;  $y_{3}\rightarrow x_{2}, x_{3}$;  $y_{1}\rightarrow x_{4}, x_{5}$;    $y_{4}\rightarrow x_{3}, x_{6}$
(b). $y_{2}\rightarrow x_{1}, x_{2}$;  $y_{3}\rightarrow x_{2}, x_{3}$;  $y_{1}\rightarrow x_{4}, x_{5}$;    $y_{4}\rightarrow x_{2}, x_{6}$
(c).  $y_{2}\rightarrow x_{1}, x_{2}$;  $y_{3}\rightarrow x_{2}, x_{3}$;  $y_{1}\rightarrow x_{4}, x_{5}$;    $y_{4}\rightarrow x_{1}, x_{6}$
(d).  $y_{2}\rightarrow x_{1}, x_{2}$;  $y_{3}\rightarrow x_{2}, x_{3}$;  $y_{1}\rightarrow x_{4}, x_{5}$;    $y_{4}\rightarrow x_{4}, x_{6}$
(e).  $y_{2}\rightarrow x_{1}, x_{2}$;  $y_{3}\rightarrow x_{3}, x_{4}$;  $y_{1}\rightarrow x_{5}, x_{6}$;    $y_{4}\rightarrow x_{1}, x_{2}$.
(f).  $y_{2}\rightarrow x_{1}, x_{2}$;  $y_{3}\rightarrow x_{3}, x_{4}$;  $y_{1}\rightarrow x_{3}, x_{6}$;    $y_{4}\rightarrow x_{1}, x_{5}$.
For  cases (a), (b) and  (c),  if $y_3$ is not adjacent to $x_4$,  then $G+y_3x_4$ contains one $K_4$, and $V(K_4)=\{y_2,y_3,x_4,x_i\}$ or  $V(K_4)=\{y_3,y_4,x_4,x_i\}$ since   $G[N_{G}(x)]$  is  $K_{3}$-free.   But  $G+xy_{2}$   or  $G+xy_{4}$   contains  two  $K_{4}$,   a contradiction.  Thus $y_{3}$  is  adjacent  to  $x_{4}$.  Similarly,  we have that $y_{3}$  is    adjacent  to  $x_{5}$.  But
$G$ contains a $F^{-}(2,4,2)$ with vertex set $\{x,x_2,x_3,x_4,x_5,y_3\}$, contradicting \autoref{lemma:2.1}.
For (d), if $y_1$ is not adjacent to $x_3$,  then $G+y_1x_3$ contains one $K_4$, and $V(K_4)=\{y_1,y_2,x_3,x_i\}$ or  $V(K_4)=\{y_1,y_4,x_3,x_i\}$  since   $G[N_{G}(x)]$  is  $K_{3}$-free.   But  $G+xy_{2}$   or  $G+xy_{4}$   contains  two  $K_{4}$,   a contradiction.  Thus $y_{1}$  is  adjacent  to  $x_{3}$.  Moreover, we can get that $y_1$  is not adjacent  to  $x_{2}$.  If not,
$G$ contains a $F^{-}(2,4,2)$ with vertex set $\{x,x_2,x_3,x_4,x_5,y_1\}$, contradicting \autoref{lemma:2.1}.
Now we claim that $y_1$  is adjacent  to  $x_{1}$.  Otherwise,  $G+y_1x_1$ contains one $K_4$ and  $x_2\notin K_4$. Moreover,  $V(K_4)=\{y_1,y_2,x_1,x_i\}$ or  $V(K_4)=\{y_1,y_4,x_1,x_i\}$   with $i\neq2$  since   $G[N_{G}(x)]$  is  $K_{3}$-free.   But  $G+xy_{2}$   or  $G+xy_{4}$   contains  two  $K_{4}$,   a contradiction.  Thus $y_{1}$  is  adjacent  to  $x_{1}$.  Similarly  as above, we can get that $y_{3}$  is  adjacent  to  $x_{5}$,  $y_3$  is not adjacent  to  $x_{4}$,  $y_{3}$  is  adjacent  to  $x_{6}$.
Moreover, we also have that   $y_{2}$  is  adjacent  to  $x_{6}$,  $y_4$  is  adjacent  to  $x_{1}$.
But  $G$ contains a $F^{-}(2,4,2)$ with vertex set $\{x_1,x_6,y_1,y_2,y_3,y_4\}$, contradicting \autoref{lemma:2.1}.
For (e), $y_{1}$  has   at most one adjacent vertex in $\{x_{1}, x_{2}\}$  since  $G$  is  $K_{4}$-free.
If  $y_{1}$  is  not adjacent  to  $x_{1}$  and  $x_{2}$,  then   $G+y_{1}x_{1}$  contains  one  $K_{4}$  and  $x_{2}\notin K_{4}$.  Moreover, we have
$V(K_4)=\{y_1,y_2,x_1,x_i\}$  or  $V(K_4)=\{y_1,y_4,x_1,x_i\}$ with $i\neq2$ in $G+y_{1}x_{1}$  since   $G[N_{G}(x)]$  is  $K_{3}$-free.  But  $G+xy_{2}$   or  $G+xy_{4}$   contains  two  $K_{4}$,   a contradiction.  Thus  $y_{1}$  has  exactly  one adjacent vertex in $\{x_{1}, x_{2}\}$.  Without loss of generality,  we assume that $y_{1}$  is adjacent to $x_{1}$.  But  $G$ contains a $W_4$ with vertex set $\{x_1,x_2,y_1,y_2,y_4\}$, contradicting \autoref{lemma:2.1}.
For (f),  if $y_2$ is not adjacent to $x_5$,  then $G+y_2x_5$ contains one $K_4$, and $V(K_4)=\{y_1,y_2,x_5,x_i\}$ or  $V(K_4)=\{y_2,y_3,x_5,x_i\}$ since   $G[N_{G}(x)]$  is  $K_{3}$-free.   But  $G+xy_{1}$   or  $G+xy_{3}$   contains  two  $K_{4}$,   a contradiction.  Thus $y_{2}$  is  adjacent  to  $x_{5}$.  But  $G+xy_2$ contains two $K_{4}$  with  vertex sets $\{x, y_2,x_1,x_2\}$ and  $\{x, y_2,x_1,x_5\}$,  a contradiction  again.

If  $\alpha=7$,  by the symmetry  of  $y_{i}$,  $i\in[4]$,  there are only  two  cases.  Without loss of generality,   we  say that they are:
(a). $y_{2}\rightarrow x_{1}, x_{2}$;  $y_{3}\rightarrow x_{2}, x_{3}$;  $y_{1}\rightarrow x_{4}, x_{5}$;    $y_{4}\rightarrow x_{6}, x_{7}$
(b). $y_{2}\rightarrow x_{1}, x_{2}$;  $y_{3}\rightarrow x_{3}, x_{4}$;  $y_{1}\rightarrow x_{6}, x_{7}$;    $y_{4}\rightarrow x_{1}, x_{5}$.
For (a), if  $y_{2}$  is  not  adjacent  to  $x_{6}$,  then  $G+y_{2}x_{6}$  contains  one  $K_{4}$,  and  $y_{1}\in K_{4}$  or  $y_{3}\in K_{4}$  since  $G[N_{G}(x)]$  is  $K_{3}$-free.  But  $G+xy_{1}$  or  $G+xy_{3}$  contains  two $K_{4}$,  a  contradiction.  Thus   $y_{2}$  is   adjacent  to  $x_{6}$.   Similarly,  we also have  that  $y_{2}$  is    adjacent  to  $x_{7}$.  But  $G$ contains a $F^{-}(2,4,2)$ with vertex set $\{x,x_1,x_2,x_6,x_7,y_2\}$, contradicting \autoref{lemma:2.1}.
For (b), if $y_2$ is not adjacent to $x_5$,  then $G+y_2x_5$ contains one $K_4$, and $V(K_4)=\{y_1,y_2,x_5,x_i\}$ or  $V(K_4)=\{y_2,y_3,x_5,x_i\}$ since   $G[N_{G}(x)]$  is  $K_{3}$-free.   But  $G+xy_{1}$   or  $G+xy_{3}$   contains  two  $K_{4}$,   a contradiction.  Thus $y_{2}$  is  adjacent  to  $x_{5}$.  But  $G+xy_2$ contains two $K_{4}$  with  vertex sets $\{x, y_2,x_1,x_2\}$ and  $\{x, y_2,x_1,x_5\}$,  a contradiction  again.

If  $\alpha=8$,   by the symmetry  of  $y_{i}$,  $i\in[4]$,  there is only  one case.  Without loss of generality,   we  say   $y_{2}\rightarrow x_{1}, x_{2}$;  $y_{3}\rightarrow x_{3}, x_{4}$;  $y_{1}\rightarrow x_{5}, x_{6}$;    $y_{4}\rightarrow x_{7}, x_{8}$.
Similarly  as  the case $\alpha=7$,  we  have that   $y_{2}$  is  adjacent  to  $x_{7}$  and  $x_{8}$.   But  $G$ contains a $F^{-}(2,4,2)$ with vertex set $\{x,x_1,x_2,x_7,x_8,y_2\}$, contradicting \autoref{lemma:2.1}.

Therefore,  $G=G_{10}$  in Figure 2   if  $e(G[A])=4$.

\textbf{Case 4.} $e(G[A])=5$.

Without loss of generality,  we  say that  $y_{1}$ is adjacent to $y_{2}$,  $y_{3}$,  $y_{4}$;  $y_{2}$ is adjacent to $y_{3}$; $y_{3}$ is adjacent to  $y_{4}$.   Moreover,  $\{y_{1}, y_{2}, y_{3}\}$ and  $\{y_{1}, y_{3}, y_{4}\}$  have  no  common adjacent vertex in  $N_{G}(x)$, respectively,   since  $G$ is  $K_{4}$-free.
$G+xy_{2}$  contains  exactly  one  $K_{4}$  and  $y_{1}, y_{3}, y_{4}\notin K_{4}$  since  $y_{1}, y_{3}, y_{4}$  are  not  adjacent to $x$.  Without loss of generality,  we  say   that  $K_{4}$ in $G+xy_{2}$  has  vertex set  $V(K_{4})=\{x,y_{2},x_{1},x_{2}\}$.  Then  $K_{4}$  in   $G+xy_{1}$  and  $K_4$  in  $G+xy_{3}$  contain at most one vertex in $\{x_{1}, x_{2}\}$  since  $G$  is  $K_{4}$-free.

If  $K_{4}$  in  $G+xy_{1}$  and $K_4$  in  $G+xy_{3}$   contain exactly  one vertex in $\{x_{1}, x_{2}\}$,  similarly as the proof  (paragraph 2) of Subcase 3.1,  it does  not hold.

If  $K_{4}$  in   $G+xy_{3}$  contains exactly  one vertex in $\{x_{1}, x_{2}\}$  and  $K_{4}$  in   $G+xy_{1}$  contains no vertex in $\{x_{1}, x_{2}\}$,
then there are two subcases. 
If  $K_{4}$  in  $G+xy_{1}$  and    $K_4$ in $G+xy_{3}$  have a common vertex  in  $N_{G}(x)$,  similarly as the proof  (paragraph 3) of Subcase 3.1,  it does  not hold.
If  $K_{4}$  in  $G+xy_{1}$  and  $K_4$  in  $G+xy_{3}$  have  no  common vertex  in  $N_{G}(x)$,  then,  without loss of generality,  let $V(K_{4})=\{x,y_{3},x_{2},x_{3}\}$  in  $G+xy_{3}$   and   $V(K_{4})=\{x,y_{1},x_{4},x_{5}\}$  in  $G+xy_{1}$.   We  have  that $y_{4}$  is  not  adjacent  to  $x_{2}$,  if  not,   $G$  contains a $W_4$ with vertex set $\{x_2,y_1,y_2,y_3,y_4\}$, contradicting \autoref{lemma:2.1}.
We claim that $y_{4}$  is   adjacent  to  $x_{3}$.  If not,   then  $G+y_{4}x_{3}$  contains  one  $K_{4}$ and  $x_{2}\notin K_{4}$.  Moreover, we have $V(K_4)=\{y_4, x_3, y_3, x_i\}$  or  $V(K_4)=\{y_4, x_3, y_1, x_i\}$ or  $V(K_4)=\{y_4, x_3, y_1, y_3\}$ with  $i\neq2$ in  $G+y_{4}x_{3}$    since $G[N_{G}(x)]$ is $K_3$-free.
If  $V(K_4)=\{y_4, x_3, y_3, x_i\}$  in  $G+y_{4}x_{3}$,  then  $G$ contains a $W_4$ with vertex set $\{x,x_2,x_3,x_i,y_3\}$, contradicting \autoref{lemma:2.1}.
If   $V(K_4)=\{y_4, x_3, y_1, x_i\}$  or  $V(K_4)=\{y_4, x_3, y_1, y_3\}$  in  $G+y_{4}x_{3}$,   then  $y_{1}$  is  adjacent to $x_{3}$.  But  $G$ contains a $W_4$ with vertex set $\{x_2,x_3,y_1,y_2,y_3\}$, contradicting \autoref{lemma:2.1} again.
Thus  $y_{4}$  is   adjacent  to  $x_{3}$.
We also claim that  $y_{4}$  is   adjacent  to  $x_{1}$.   If  not,  then  $G+y_{4}x_{1}$  contains  one  $K_{4}$   and  $x_{2}\notin K_{4}$.  Moreover, we have $V(K_4)=\{y_4, x_1, y_3, x_i\}$  or  $V(K_4)=\{y_4, x_1, y_1, x_i\}$ or  $V(K_4)=\{y_4, x_1, y_1, y_3\}$ with  $i\neq2$ in  $G+y_{4}x_{1}$    since $G[N_{G}(x)]$ is $K_3$-free.    If   $V(K_4)=\{y_4, x_1, y_3, x_i\}$ or  $V(K_4)=\{y_4, x_1, y_1, y_3\}$  in  $G+y_{4}x_{1}$,  then  $y_{3}$  is  adjacent to $x_{1}$.  But  $G$  contains a   $K_{4}$  with  vertex set  $\{y_{2},y_{3},x_{1},x_{2}\}$, a contradiction.  If  $V(K_4)=\{y_4, x_1, y_1, x_i\}$  in  $G+y_{4}x_{1}$,  then  $G+xy_{1}$  contains two $K_{4}$  with  vertex sets $\{x, y_1,x_4,x_5\}$ and  $\{x, y_1,x_1,x_i\}$,  a contradiction  again. Thus  $y_{4}$  is   adjacent  to  $x_{1}$.
Since  $G$  is  $K_{4}$-free,  $y_{2}$  and  $y_{4}$  have at most one  adjacent vertex  in $\{x_{4}, x_{5}\}$.
By a similar reasoning to the proof (paragraph 4) of  Subcase 3.1,   $y_{2}$  and  $y_{4}$   have  exactly  one  adjacent vertex in  $\{x_{4}, x_{5}\}$. Moreover, we have that  $y_{2}$  and  $y_{4}$   have different  adjacent vertex in  $\{x_{4}, x_{5}\}$  since  $G+y_{2}y_{4}$  contains  exactly one  $K_{4}$.   Without loss of generality,  let  $y_{2}$  be adjacent to $x_{5}$,  $y_{4}$  be adjacent to $x_{4}$.  And  $y_{1}$  is  not adjacent  to  $x_{1}$.  Otherwise,  $G$ contains a $W_4$ with vertex set $\{x_1,y_1,y_2,y_3,y_4\}$, contradicting \autoref{lemma:2.1}.
By the above discussion  and  $G$  is  $K_{4}$-free,  we can get that  for  vertex set $\{x_{i}:i\in[5]\}$,   $y_{1}$  is only adjacent to  $x_{4}$,  $x_{5}$;   $y_{2}$  is only adjacent to  $x_{1}$, $x_{2}$,  $x_{5}$;  $y_{3}$  is only adjacent to  $x_{2}$, $x_{3}$;  $y_{4}$  is only adjacent to  $x_{1}$, $x_{3}$,  $x_{4}$.

Since $y_{1}$  is  not adjacent  to  $x_{1}$,   $G+y_{1}x_{1}$  contains  one  $K_{4}$  and  $x_{2}, y_3\notin K_{4}$.  Moreover,   we have
$V(K_4)=\{y_1,x_1,y_2,x_i\}$  or  $V(K_4)=\{y_1,x_1,y_4,x_i\}$ with $i\neq2$ in $G+y_{1}x_{1}$  since  $G[N_{G}(x)]$ is $K_3$-free  and  $G$ is $K_4$-free.
If  $V(K_4)=\{y_1,x_1,y_2,x_i\}$  in  $G+y_{1}x_{1}$,   then  $G$ contains a $W_4$ with vertex set $\{x,x_1,x_2,x_i,y_2\}$, contradicting \autoref{lemma:2.1}.
%$G+xy_{2}$  has two $K_{4}$ with vertex sets $\{x,y_2,x_1,x_2\}$  and  $\{x,y_2,x_1,x_i\}$,  a contradiction.
Thus    $G+y_{1}x_{1}$  contains  one  $K_{4}$  with vertex set $V(K_4)=\{y_1,x_1,y_4,x_i\}$  and  $i\neq2, 3, 5$.    Now  we claim $i=4$.  Otherwise,  without loss of generality,  we  say   that  $K_{4}$  in  $G+y_{1}x_{1}$  has  vertex set  $\{y_{1}, x_{1},y_{4},x_{6}\}$.  That is,   $x_6$  is adjacent to $y_1$  and $y_4$.
Thus  $x_{6}$  is  not adjacent to $y_{2}$  since  $G+y_{2}y_{4}$  contains exactly  one  $K_{4}$.  Then  $G+y_{2}x_{6}$  contains one  $K_{4}$ and  $y_3\notin K_4$ since  $G$ is $K_4$-free.  Moreover, we have $V(K_4)=\{y_2,x_6,y_1,x_j\}$  in $G+y_{2}x_{6}$ since $G[N_{G}(x)]$  is $K_3$-free.  But $G+xy_{1}$  contains  two  $K_{4}$ with vertex sets $\{x,y_1,x_4,x_5\}$  and  $\{x,y_1,x_6,x_j\}$,  a contradiction.
Therefore, $i=4$, that is,   $K_{4}$  in  $G+y_{1}x_{1}$  has  vertex set  $\{y_{1}, x_{1},y_{4},x_{4}\}$,  thus  $x_{1}$  is adjacent to  $x_{4}$.
If  $x_{3}$  is adjacent  to  $x_{4}$,  then  $G$ contains a $W_4$ with vertex set $\{x,x_1,x_2,x_3,x_4\}$, contradicting \autoref{lemma:2.1}.
Thus   $x_{3}$  is  not  adjacent  to  $x_{4}$  and   $G+x_{3}x_{4}$  contains  one  $K_{4}$.  Since  $x_3$ and $x_4$ have only one common adjacent vertex $y_4$  in $\{y_i, i\in[4]\}$,  we can get $y_4\notin K_4$  in $G+x_{3}x_{4}$.  If not,   $G+xy_4$ contains  two  $K_4$, a contradiction.  Thus  $K_{4}$  in  $G+x_{3}x_{4}$  has  vertex set  $\{x, x_{3},x_{4},x_{p}\}$   with  $p\neq1, 2$.  Now  we claim $p=5$.  Otherwise,  without loss of generality,  we  say   that
$K_{4}$  in  $G+x_{3}x_{4}$  has  vertex set  $\{x, x_{3},x_{4},x_{6}\}$.  That is,  $x_6$  is adjacent to $x_3$  and $x_4$. Thus for  $i=1,3,4$,  $y_{i}$  is  not adjacent to  $x_{6}$   since   $G+xy_{i}$  contains exactly  one  $K_{4}$.   But  $G+y_{4}x_{6}$  contains  one $K_{4}$  with vertex set $\{y_4,x_6,x_s,x_t\}$,  which contradicts  the fact that $G[N_{G}(x)]$ is $K_3$-free.  Therefore, $p=5$, that is,  $K_{4}$  in  $G+x_{3}x_{4}$  has  vertex set  $\{x, x_{3},x_{4},x_{5}\}$  and  $x_{3}$  is adjacent to  $x_{5}$.
By the above discussion  and $G$  is  $K_{4}$-free,  we can get that  $G[\{x_i, i\in[5]\}]=C_5$  with  $C_5=x_1x_2x_3x_5x_4x_1$.

If  $|N_{G}(x)|=5$,  which  means $N_{G}(x)=\{x_{i}: i\in[5]\}$,  then  $G=G_{10}$  in Figure 2 with $a_0=x_1, a_1=y_4,  a_2=x_3, a_3=x_5,  a_4=y_2,  b_0=x,  b_1=x_4,  b_2=y_1,  b_3=y_3,  b_4=x_2$  by  a routine check.

If  $|N_{G}(x)|\geq6$, then  $G[\{x, y_{i}, x_{j}: i\in[4], j\in[5]\}]=G_{10}$  in  Figure 2.
If  $x_{4}$  is adjacent to  $x_{6}$,  then  $y_{1}$  and  $y_{4}$  are not adjacent to  $x_{6}$  since  $G+xy_{1}$  and  $G+xy_{4}$  contain exactly   one  $K_{4}$.
For the $K_4$  in $G+y_{4}x_{6}$,  we have $y_{1}\notin K_4$ since $y_{1}$  is not adjacent to  $x_{6}$. Then
$y_{3}\in K_4$ since $G[N_{G}(x)]$  is $K_3$-free.  That is,  $K_4$ in $G+y_{4}x_{6}$  has vertex set $\{y_4,x_6,y_3,x_i\}$.  But  $G+xy_3$  contains  two $K_4$ with vertex sets $\{x,y_3,x_2,x_3\}$  and $\{x,y_3,x_6,x_i\}$,  a contradiction.  Thus $x_{4}$  is  not adjacent to  $x_{6}$.  Similarly,   we  have  that  $x_{2}$, $x_{3}$, $x_{5}$  are  not adjacent to  $x_{6}$.

If  $y_{2}$  is  not  adjacent  to  $x_{6}$,  then  $G+y_{2}x_{6}$  contains  one $K_{4}$  and  $\{y_{1}, y_{3}\}\subset K_{4}$   since $G[N_{G}(x)]$ is $K_3$-free  and  $G+xy_{1}$,  $G+xy_{3}$  contain  exactly  one  $K_{4}$.    Thus  $y_{1}$  and  $y_{3}$  are adjacent  to  $x_{6}$,  $y_{4}$  is  not adjacent  to  $x_{6}$ since $G$  is $K_4$-free.
For  the $K_4$ in $G+x_{3}x_{6}$,  since  $x_3$  and $x_6$  have only one common adjacent vertex $y_3$  in $\{y_i, i\in[4]\}$,  we can get $y_3\notin K_4$  in $G+x_{3}x_{6}$.  If not,   $G+xy_3$  contains two  $K_4$, a contradiction.  Thus  $K_{4}$  in  $G+x_{3}x_{6}$  has  vertex set  $\{x, x_{3},x_{6},x_{p}\}$   with  $p\geq7$.   Without loss of generality,  let  $V(K_4)=\{x, x_{3},x_{6},x_{7}\}$  in  $G+x_{3}x_{6}$,  then $x_{3}$  and  $x_{6}$  are  adjacent to  $x_{7}$.
Thus for  $i=1,3,4$,  $y_{i}$  is  not adjacent to  $x_{7}$   since   $G+xy_{i}$  contains exactly  one  $K_{4}$.   But  $G+y_{4}x_{7}$  contains  one $K_{4}$  with vertex set $\{y_4,x_7,x_s,x_t\}$,  which contradicts  the fact that $G[N_{G}(x)]$ is $K_3$-free.

Now  we assume that   $y_{2}$  is   adjacent  to  $x_{6}$.
If  $y_{4}$  is  not   adjacent  to  $x_{6}$,   then  $G+y_{4}x_{6}$  contains  one $K_{4}$  and  $\{y_{1}, y_{3}\}\subset K_{4}$   since  $G[N_{G}(x)]$ is $K_3$-free  and  $G+xy_{1}$,  $G+xy_{3}$  contain  exactly  one  $K_{4}$.  Thus  $y_{1}$  and  $y_{3}$  are adjacent  to  $x_{6}$.    But  $G$  contains  $K_{4}$  with  vertex set  $\{y_{1}, y_{2}, y_{3}, x_{6}\}$, a contradiction. Thus  $y_{4}$  is   adjacent  to  $x_{6}$.
Moreover,  $y_{1}$  and  $y_{3}$  are  not adjacent to  $x_{6}$  since  $G+y_{2}y_{4}$ contains exactly one copy of $K_{4}$.
For  the $K_4$ in $G+x_{3}x_{6}$,  since  $x_3$  and $x_6$ have only one common adjacent vertex $y_4$  in $\{y_i, i\in[4]\}$,  we can get $y_4\notin K_4$  in $G+x_{3}x_{6}$.  If not,   $G+xy_4$   contains  two  $K_4$, a contradiction.  Thus  $K_{4}$  in  $G+x_{3}x_{6}$  has  vertex set  $\{x, x_{3},x_{6},x_{p}\}$   with  $p\geq7$.   Without loss of generality,  let  $V(K_4)=\{x, x_{3},x_{6},x_{7}\}$  in  $G+x_{3}x_{6}$,  then $x_{3}$  and  $x_{6}$  are  adjacent to  $x_{7}$.
Thus for  $i=2,3,4$,  $y_{i}$  is  not adjacent to  $x_{7}$   since   $G+xy_{i}$  contains  exactly  one  $K_{4}$.
For the $K_4$ in $G+y_{4}x_{7}$,
we have $y_{3}\notin K_4$ since $y_{3}$  is not adjacent to  $x_{7}$.   Then
$y_{1}\in K_4$  in $G+y_{4}x_{7}$  since $G[N_{G}(x)]$   is  $K_3$-free.   That is,  $K_4$ in $G+y_{4}x_{7}$ has vertex set $\{y_4,x_7,y_1,x_i\}$.  But  $G+xy_1$  contains  two $K_4$ with vertex sets $\{x,y_1,x_4,x_5\}$  and $\{x,y_1,x_7,x_i\}$,  a contradiction.

Therefore,  if   $K_{4}$  in   $G+xy_{3}$  contains exactly  one vertex in $\{x_{1}, x_{2}\}$  and  $K_{4}$  in   $G+xy_{1}$  contains no vertex in $\{x_{1}, x_{2}\}$,  then   we can get $G=G_{10}$  in Figure 2.

If  $K_{4}$  in   $G+xy_{1}$  contains exactly  one vertex in $\{x_{1}, x_{2}\}$  and  $K_{4}$  in   $G+xy_{3}$  contains no vertex in $\{x_{1}, x_{2}\}$,
then, by the symmetry of $y_1$  and  $y_3$,  we can get $G=G_{10}$  in Figure 2.

If $K_{4}$  in  $G+xy_{1}$ and  $K_4$   in  $G+xy_{3}$  contain  no vertex in $\{x_{1}, x_{2}\}$,
then there are two subcases.
If  $K_{4}$  in  $G+xy_{1}$  and  $K_4$  in   $G+xy_{3}$  have a common vertex  in  $N_{G}(x)$,   similarly as the proof   (paragraph 5) of Subcase 3.1,  it does  not hold.
If  $K_{4}$  in  $G+xy_{1}$  and  $K_4$  in $G+xy_{3}$  have no  common vertex  in  $N_{G}(x)$,  then,  without loss of generality,  let  $V(K_{4})=\{x,y_{3},x_{3},x_{4}\}$  in  $G+xy_{3}$   and   $V(K_{4})=\{x,y_{1},x_{5},x_{6}\}$  in  $G+xy_{1}$.
If  $y_{4}$  is  not adjacent to  $x_{1}$,  then  $G+y_{4}x_{1}$  contains one $K_{4}$.   Moreover,  we have $V(K_4)=\{y_4, x_1, y_3, x_i\}$  or  $V(K_4)=\{y_4, x_1, y_1, x_i\}$ or  $V(K_4)=\{y_4, x_1, y_1, y_3\}$  in  $G+y_{4}x_{1}$    since $G[N_{G}(x)]$ is $K_3$-free.
If   $V(K_4)=\{y_4, x_1, y_1, x_i\}$ or  $V(K_4)=\{y_4, x_1, y_3, x_i\}$  in  $G+y_{4}x_{1}$,  then $G+xy_{1}$   or  $G+xy_{3}$  contains  two  $K_{4}$,  a contradiction.  If  $V(K_4)=\{y_4, x_1, y_1, y_3\}$  in  $G+y_{4}x_{1}$,  then $y_{1}$  and  $y_{3}$  are  adjacent  to  $x_{1}$.    But  $G$  contains   $K_{4}$  with  vertex set  $\{y_{1},y_{2},y_{3},x_{1}\}$, a contradiction  again.
Thus  $y_{4}$  is  adjacent to $x_{1}$.  Similarly,  we can get that $y_{4}$  is also adjacent to $x_{2}$.
But  $G$ contains a $F^{-}(2,4,2)$   with vertex set $\{x_1,x_2,y_1,y_2,y_3,y_4\}$, contradicting \autoref{lemma:2.1}.

Therefore,  $G=G_{10}$  in Figure 2   if  $e(G[A])=5$.

In conclusion,  $G_{10}$  in Figure 2  is the only  uniquely $K_4$-saturated graph  on  $n$  vertices  with  maximum degree  $n-5$.
\end{proof}

\section*{Declaration of competing interest}

The authors declare that they have no known competing financial interests or personal relationships that could have
appeared to influence the work reported in this paper.

\section*{Acknowledgements}

The work was supported by  Research Foundation of China Scholarship Council (No.\ 202608320356), National Natural Science Foundation of China (No.\ 12271251), Postgraduate Research \& Practice Innovation Program of Jiangsu Province (No. 26CXJH0986) and Fundamental Research Funds for the Central
Universities (Grant No. NJ2026048).

\section*{Data availability}

No data was used for the research described in the article.


\begin{thebibliography}{99}


\bibitem{L.W.Berman2016}   L.W. Berman, G.G. Chappell, J.R. Faudree, J. Gimbel,  C. Hartman,
Uniquely tree-saturated graphs, Graphs Combin.  \textbf{32}   (2016)  463--494.

\bibitem{A.Bondy2008} J.A. Bondy,   U.S.R. Murty,  Graph Theory,  Springer,  London,  2008.

\bibitem{J.Cooper2012}  J. Cooper,   J. Lenz,  T.D.  LeSaulnier,   P.S.  Wenger,  D.B.  West,
Uniquely  $C_{4}$-saturated graphs, Graphs Combin.  \textbf{28}  (2012)  189--197.

\bibitem{B.L.Currie2021} B.L. Currie,  J.R. Faudree,  R.J. Faudree,  J.R. Schmitt,
A survey of minimum saturated graphs, Electron. J. Comb.   (2021)  \#DS19.


\bibitem{P.Erdos1964} P. Erd\H{o}s, A. Hajnal, J.W. Moon,  A problem in graph theory,
Amer. Math. Monthly  \textbf{71}  (1964) 1107--1110.


\bibitem{A.Gyarfas2018} A. Gy\'{a}rf\'{a}s,  S.G. Hartke,   C. Viss,   Uniquely $K_{r}^{(k)}$-saturated hypergraphs,
Electron. J. Comb.  \textbf{25}  (2018)  \#P4.35.



\bibitem{S.G.Hartke2012}  S.G.  Hartke,  D.  Stolee,  Uniquely $K_{r}$-saturated graphs,
Electron. J. Comb.  \textbf{19}  (2012) \#P6.

\bibitem{Huang2024} S. Huang, H. Lei, Y. Shi,  J. Zhang,
The saturation number of $K_{3,3}$,
Discrete Math. \textbf{347} (2024) 113794.


\bibitem{Ji2026} X. Ji, E.R. van Dam,  J. Guo, Z. Wang, Bowtie-saturated graphs and the spectral radius,
Discrete Math.
\textbf{349} (2026) 115254.

\bibitem{Lan2025} Y. Lan, Y. Shi, Y. Wang,  J. Zhang,
The saturation number of $C_6$,
Discrete Math. \textbf{348} (2025) 114504.


 \bibitem{LOSWZ2019} H. Lei, S. O, Y. Shi, D.B. West, X. Zhu,
Extremal problems on saturation for the family of $k$-edge-connected graphs,
 Discrete Appl. Math. \textbf{260} (2019) 278--283.



\bibitem{Li2026}  Y. Li, K. Xu, D. Gerbner, W. Liu,
Uniquely $C_{4}^{+}$-saturated graphs,  Discrete Math.  \textbf{349} (2026) 114874.



\bibitem{W.Mantel1907}  W. Mantel, Problem 28, soln. by H. Gouwentak, W. Mantel, J. Teixeira de Mattes, F. Schuh and W.A. Wythoff,
Wiskundige Opgaven  \textbf{10}  (1907) 60--61.

\bibitem{P.Turan1941} P. Tur\'{a}n,
Eine Extremalaufgabe aus der Graphentheorie,
Mat. Fiz. Lapok \textbf{48}   (1941) 436--452.


\bibitem{P.S.Wenger2016}    P.S. Wenger,  D.B. West,
Uniquely cycle-saturated graphs,
J. Graph Theory  \textbf{85}  (2017) 94--106.


\bibitem{P.S.Wenger2010}
P.S. Wenger, Three existence problems in extremal graph theory,
Ph.D. Thesis, University of Illinois at Urbana-Champaign, 2010.

\end{thebibliography}
\end{document}